\documentclass[12pt]{article}
\usepackage{geometry}
\usepackage{amsmath,amsthm,amscd,amssymb,latexsym,upref,enumerate,bbm,bm,lscape,epsfig,nicefrac,hyperref,amsfonts, color,graphicx,xcolor,cite,mathtools}

\usepackage{comment}

\newtheorem{info}{}
\newtheorem{theorem}[info]{Theorem}
\newtheorem{lemma}[info]{Lemma}
\newtheorem{proposition}[info]{Proposition}
\newtheorem{corollary}[info]{Corollary}

\newtheorem{remark}[info]{Remark}

\numberwithin{info}{section}
\renewcommand{\[}{\begin{equation}}
	\renewcommand{\]}{\end{equation}}
\makeatletter
\g@addto@macro\normalsize{%
	\setlength\abovedisplayskip{5pt}
	\setlength\belowdisplayskip{5pt}
	\setlength\abovedisplayshortskip{4pt}
	\setlength\belowdisplayshortskip{4pt}}
\makeatother

\newcommand{\N}{\mathbb{N}}
\newcommand{\C}{\mathbb{C}}
\newcommand{\R}{\mathbb{R}}

\newcommand{\E}{\mathbb{E}}
\renewcommand{\P}{\mathbb{P}}
\newcommand{\bH}{\mathbb{H}}

\newcommand{\lam}{\lambda}

\renewcommand{\cal}{\mathcal}

\newcommand{\Tr}{\mathrm{Tr}}

\newcommand{\sq}{\sqrt}

\renewcommand{\sq}{\sqrt}

\newcommand{\Sum}{\mathrm{\Sum}}

\renewcommand{\Im}{\mathrm{Im}}

\newcommand{\Crit}{\mathrm{Crit}}

\newcommand{\D}{\nabla}
\newcommand{\vol}{\mathrm{vol}}

\newcommand{\n}{\textbf{n}}
\newcommand{\W}{\mathrm{W}}
\newcommand{\GOE}{\mathrm{GOE}}
\newcommand{\Cov}{\mathrm{Cov}}

\newcommand{\sign}{\mathrm{sign}}

\title{The Ground State Energy and Concentration of Complexity in Spherical Bipartite Models}

\author{
	Pax Kivimae\thanks{Department of Mathematics, Northwestern University. Email: kivimae@math.northwestern.edu, }
}

\date{\today}

\begin{document}
	
	\maketitle
	
	\begin{abstract}
		We establish an asymptotic formula for the ground-state energy of the spherical pure $(p,q)$-spin glass model for $p,q\ge 97$. We achieve this by showing a concentration result for the complexity of critical points with values within a region of the ground state energy. More specifically, we show that the second moment of this count coincides with the square of the first moment up to a sub-exponential factor.
	\end{abstract}

	\section{Introduction}

	In recent years, bipartite and, more generally, multi-species spin models have seen increased interest in both mathematical and physical literature, owing in part to their applications to the theory of neural networks \cite{cs1,cs2,cs3}, and to theoretical biology, \cite{bio1,bio2,bio3,bio4}, among others \cite{socio1,socio2}. This work will consider the spherical pure $(p,q)$-spin bipartite spin glass model introduced in \cite{tuca}. The Hamiltonian of this model $H_{p,q,N_1,N_2}:S^{N_1-1}(\sq{N_1})\times S^{N_2-1}(\sq{N_2})\to \R$ is defined for $(\sigma,\tau)\in S^{N_1-1}(\sq{N_1})\times S^{N_2-1}(\sq{N_2})$ by
	\[H_{p,q,N_1,N_2}(\sigma,\tau)=\sum_{i_1,\dots i_p=1}^{N_1}\sum_{j_1,\dots j_q=1}^{N_2}J_{i_1,\dots i_p;j_1,\dots j_q}\sigma_{i_1}\dots \sigma_{i_p}\tau_{j_1}\dots \tau_{j_q},\]
	where here $J_{i_1,\dots i_p;j_1,\dots j_q}$ are i.i.d centered Gaussian random variables with variance $N/(N_1^pN_2^q)$ and $S^{n}(m)=\{\sigma\in \R^n:\|\sigma\|=m\}$. Equivalently this may be defined as the unique smooth centered Gaussian field on $S^{N_1-1}(\sq{N_1})\times S^{N_2-1}(\sq{N_2})$ with covariance function defined for $(\sigma,\tau),(\sigma',\tau')\in S^{N_1-1}(\sq{N_1})\times S^{N_2-1}(\sq{N_2})$ by
	\[\E[ H_{p,q,N_1,N_2}(\sigma,\tau)H_{p,q,N_1,N_2}(\sigma',\tau')]=N(\frac{1}{N_1}(\sigma,\sigma'))^p(\frac{1}{N_2}(\tau,\tau'))^q,\]
	where on the right-hand side $(\sigma,\sigma')$ and $(\tau,\tau')$ denote the Euclidean inner product. For any choice of $0<\gamma<1$, we further consider a choice, for each $N$, of $N_1,N_2\ge 2$, such that $N_1+N_2=N$ and such that $|\gamma N-N_1|=O(1)$. With such a sequence fixed, we will denote $H_{N}:=H_{p,q,N_1,N_2}$. 
	
	The asymptotics for the annealed complexity of $H_N$ have recently been obtained in \cite{bipartite}, whose results we now recall. For a Borel subset $B\subseteq \R$, let us denote by $\Crit_{N}(B)$ the count of critical points of $H_N$ such that $ H_N\in NB$. We let $\Crit_{N,0}(B)$ denote the count of the subset of such critical points that are also local minima. Now the results of \cite{bipartite} show that there is a fixed function $\Sigma_{p,q,\gamma}$ (recalled below), and a fixed value $E_{p,q,\gamma;\infty}>0$, such that for any $E\in \R$, we have that
	\[\lim_{N\to \infty}\frac{1}{N}\log(\E[\Crit_N((-\infty,E))])=\sup_{t<E}\Sigma_{p,q,\gamma}(t),\label{Single point complexity}\]
	\[\lim_{N\to \infty}\frac{1}{N}\log(\E[\Crit_{N,0}((-\infty,E))])=\sup_{t<\min(E,-E_{p,q,\gamma;\infty})}\Sigma_{p,q,\gamma}(t).\label{Single point complexity minima}\]
	From Theorem 2 of \cite{tuca}, one may conclude that $\Sigma_{p,q,\gamma}(t)$ is somewhere positive on 
	\linebreak
	$(-\infty,-E_{p,q,\gamma;\infty})$, and is negative for sufficiently small $t$. Let us denote by $E_{p,q,\gamma;0}$, the largest value such that $\Sigma_{p,q,\gamma}(t)<0$ for $t<-E_{p,q,\gamma;0}$. We now state our main result.
	\begin{theorem}
		\label{Main-Theorem}
		Let us assume that $p,q\ge 97$ and $0<\gamma<1$. Then for any choice of $E\in (-E_{p,q,\gamma;0},-E_{p,q,\gamma;\infty})$ we have that
		\[\lim_{N\to\infty}\frac{1}{N}\log(\frac{\E[\Crit_N((-\infty,E))^2]}{\E[\Crit_N((-\infty,E))]^2})=0.\label{main-theorem-eqn}\]
	\end{theorem}
	
	Similar results have been obtained for the unipartite spherical pure $p$-spin glass model in \cite{pspin-second} and for certain mixed spherical models in \cite{subag-1RSB}. We also mention the works \cite{nic1,nic2,second1,second2}, in which the concentration of complexity is understood for a variety of Gaussian fields through similar methods. 
	
	Now, using Theorem \ref{Main-Theorem}, and an application of the Paley-Zigmund inequality similarly to that of Appendix IV of \cite{pspin-second}, we may conclude the following result.
	\begin{corollary}
		For $p,q\ge 97$ and $0<\gamma<1$, we a.s. have that
		\[\lim_{N\to \infty} \frac{1}{N}\min_{(\sigma,\tau)}H(\sigma,\tau)=-E_{p,q,\gamma;0}.\]
	\end{corollary}
	
	We now comment on some related results in spin glass theory. The ground-state energy, and more generally the free energy, has been rigorously computed for mixed Ising and spherical spin glasses \cite{Parisi-Even-Ising, Parisi-Even-Spherical, Parisi-General-Ising, Parisi-General-Spherical}, vector spin glasses \cite{Parisi-Potts-Spin, Parisi-Vector-Spin}, and a variety of others. In the case of the bipartite spin glass model, though, the free energy has only been computed for models in which the interaction strength is positive definite \cite{Bipartite-PD}, at high temperature \cite{tuca}, or in the case of the pure $(1,1)$-spin model for non-critical temperatures \cite{11model}.
	
	A key problem in extending these computations to bipartite models is the inability to apply the Guerra interpolation scheme \cite{Guerra}, which plays a key role in the derivations in the unipartite case. Our method sidesteps this hurdle by instead relying wholly on complexity-based methods. In particular, the methods of this paper are in part motivated by those of \cite{pspin-second,subag-1RSB}, which establish concentration results on the complexity of both pure spherical spin glasses and a subsector of mixed spherical spin glasses within the $1$-RSB phase. In particular, they are also able to rederive the Crisanti-Sommers formula for the ground-state energy without the use of the Guerra interpolation scheme.
	
	We also mention that as in the derivation of the first moment of complexity in \cite{bipartite}, we rely crucially on recent results on the matrix Dyson equation \cite{MDE1,MDE2}, which allows us to understand the limit of the average spectral measure of certain Gaussian random matrices appearing as Hessians in our analysis.
	
	The method of proof is structurally similar to the approach developed in \cite{pspin-second} to treat the unipartite model. The primary differences will involve the analysis of certain deterministic functions which appear in the course of our approach. Often, we find that these functions are sufficiently complicated that we are unable to adapt the direct analytic proofs of \cite{pspin-second}, which work in the unipartite case, as well as the inability to rely on numerical methods to treat marginal cases here (note, for example, that strengthening a bound numerically from $p,q\le 97$ to $p,q\le 96$ involves the verification of an infinite number of cases, even ignoring the continuous values taken by $\gamma$). Foremost among these issues is our limited ability to understand the limiting empirical spectral measure of the Hessian. Indeed, while one is able to implicitly determine this measure through the methods of \cite{MDE1}, even basic properties, such as the value of the largest point in support of the measure, are relatively unclear, in comparison to the exact computations which are available in the case of the semicircle law which appears in the unipartite case \cite{pspin-one}. 
	
	The main new tools we need to employ in our analysis are relations that we obtain by analyzing the relevant quantities for finite $N$. Indeed, by rescaling the entries of the Hessian in a certain fixed way, we are able to change the determinant of the Hessian by a fixed multiple while simplifying the energy dependence of the limiting empirical spectral measures given in \cite{bipartite} so that all measures are simple translates of the same measure. This not only allows us to remove multiple variables from our minimization, but it allows us to effectively compare this determinant to the semi-circle law by using classical matrix inequalities. In particular, both of these methods result in inequalities involving fixed limiting functions but which appear relatively opaque from this perspective.
	
	We also note that we believe that Theorem \ref{Main-Theorem} holds under the weaker restriction of $p,q\ge 2$. The requirement that $p,q\ge 97$ results from difficulties in controlling fixed functions occurring in the complexity. In particular, it is clear from the results of Section \ref{section:outline}, that for fixed $p,q\ge 5$ and $0<\gamma<1$, one may obtain (\ref{main-theorem-eqn}) by showing that for $E\in (-E_{p,q,\gamma;0},-E_{p,q,\gamma;\infty})$, we have that \[\sup_{r,t\in (-1,1),u,v\le E}\Sigma_{p,q,\gamma,2}(r,t,u,v)=2\Sigma_{p,q,\gamma}(E),\]
	where $\Sigma_{p,q,\gamma,2}$ is a fixed function introduced in (\ref{Sigma_2:def}). In particular, given a fixed choice of $(p,q,\gamma)$, one may check this equality using numerical methods. On the other hand, even when one fixes a value of $(p,q)$, as the values of $\gamma$ are continuous, it is not clear how to incorporate a single additional value of $(p,q)$ as it is unclear how to obtain an explicit expression for the limiting empirical spectral density which are computationally manageable. Indeed, improvement of the theorem to this point would likely require explicit information about the limiting empirical measure, a result that does not appear technically feasible outside of the semi-circle case.
	
	With these technicalities noted, we now describe the structure of the paper. In Section \ref{section:outline} we introduce a sequence of lemmas from which we prove Theorem \ref{Main-Theorem}. Sections \ref{section:lemma-section} and \ref{section:KR-2point} are devoted to the proofs of these lemmas. In addition, we will provide three appendices. In Appendix \ref{section:covariance-section} we will provide a number of computations related to the covariance structure of $H_N$, and its derivatives, at two points. In Appendix \ref{section:appendix-KR} we will provide a proof of Proposition \ref{prop:KR-2point}, which is our application of Kac-Rice formula to the computation of $\E[\Crit_N(B)^2]$. Lastly, in Appendix \ref{appendix:mckenna rescaling} we will relate our definition of $\Sigma_{p,q,\gamma}$ provided below to the one present in \cite{bipartite}.
	
	\subsection{Notation\label{subsection:notation}}
	We define for a Borel subset $B\subseteq \R$,
	\[\Crit_N(B)=\#\{z\in S^{N_1-1}(\sq{N_1})\times S^{N_2-1}(\sq{N_2}):\D H_N(z)=0,H_N(z)\in NB\}.\]
	Here $\D$ denotes the gradient taken with respect to the standard Riemannian metric on the sphere. We define for Borel subsets $B\subseteq \R$ and $I_O\subseteq [-1,1]^2$, 
	\[\Crit_{N,2}(B,I_O)=\#\{(\sigma,\tau),(\sigma',\tau')\in (S^{N_1-1}(\sq{N_1})\times S^{N_2-1}(\sq{N_2}))^2: \]\[\D H_N(\sigma,\tau)=\D H_N(\sigma',\tau')=0, H_N(\sigma,\tau),H_N(\sigma',\tau')\in NB,\;\;(\frac{1}{N_1}(\sigma,\sigma'),\frac{1}{N_2}(\tau,\tau'))\in I_O\}.\]
	We observe that $\Crit_{N,2}(B,[-1,1]^2)=\Crit_N(B)^2$. We will say that a subset $B\subset \R$ is nice if it is a finite union of open intervals, and similarly, we will say that a subset $I_O\subset [-1,1]^2$ is nice if it is a finite union of products of open intervals.
	
	We now recall the functions composing the definition of $\Sigma_{p,q,\gamma}$. We will denote by $\bH=\{z\in \C:\Im(z)>0\}$ the complex upper half plane. Let $(m_0(z),m_1(z))$ denote the unique functions $m_i:\bH\to \bH$ satisfying the following system of equations
	\[\begin{cases}
		1+(z+\gamma p^{-1}(p-1)m_0(z)+(1-\gamma)m_1(z))m_0(z)=0\\
		1+(z+\gamma m_0(z)+(1-\gamma) q^{-1}(q-1)m_1(z))m_1(z)=0\end{cases}.\label{system of equations}\]
	There are unique compactly supported measures, $\mu_{p,q,\gamma,i}$, whose Stiejes transform at $z$ is given by $m_i(z)$, for $i=0,1$. Given these we define the measure $\mu_{p,q,\gamma}=\gamma \mu_{p,q,,\gamma,0}+(1-\gamma)\mu_{p,q,\gamma,1}$. We now define
	\[\Omega_{p,q,\gamma }(E)=\int \log(|x-E|)\mu_{p,q,\gamma}(dx),\;\; \Theta_{p,q,\gamma}(E)=-\frac{1}{2}E^2+\Omega_{p,q,\gamma}(E),\label{eqn:omega-def}\]
	\[C_{p,q,\gamma}=\frac{1}{2}(1+\gamma\log(\frac{p}{\gamma})+(1-\gamma)\log(\frac{q}{(1-\gamma)})),\;\;\Sigma_{p,q,\gamma}(E)=C_{p,q,\gamma }+\Theta_{p,q,\gamma}(E).\label{eqn:c-def}\]
	\begin{remark}
		The existence of such $\mu_{p,q,\gamma}$ and the verification that our $\Sigma_{p,q,\gamma}$ coincides with the corresponding quantity in \cite{bipartite} are given in Appendix \ref{appendix:mckenna rescaling}.
	\end{remark}
	
	We now define some functions that will show up in the computation of  $\E[\Crit_{N,2}(B,I)]$. For $r,t\in (-1,1)$, we define
	\[g_{p,q,\gamma}(r,t)=\frac{\gamma }{2}\log(\frac{1-r^2}{1-t^{2q}r^{2p-2}})+\frac{(1-\gamma)}{2}\log(\frac{1-t^2}{1-t^{2q-2}r^{2p}}),\label{h-def}\]
	\[\Psi_{p,q,\gamma}(r,t,E_1,E_2)=g_{p,q,\gamma}(r,t)-\frac{1}{2}(E_1,E_2)\Sigma_U(r,t)^{-1}(E_1,E_2)^t+\Omega_{p,q,\gamma}(E_1)+\Omega_{p,q,\gamma}(E_2),\label{eqn:psi-def}\]
	where $\Sigma_U(r,t)$ is a certain invertible $2$-by-$2$ matrix introduced in (\ref{appendix Sigma_E}) (see also (\ref{better U-def})). We further define
	\[\Sigma_{p,q,\gamma,2}(r,t,E_1,E_2)=2C_{p,q,\gamma}+\Psi_{p,q,\gamma}(r,t,E_1,E_2).\label{Sigma_2:def}\]
	We define a special function that will appear in our analysis. \[\Omega_{sc}(E)=\int \log(|x-E|)\frac{\sq{4-x^2}}{2\pi}dx=\]
	\[\frac{E^2}{4}-\frac{1}{2}-\begin{cases}\frac{|E|}{4}\sq{E^2-4}-\log(\frac{|E|+\sq{E^2-4}}{2});|E|\ge 2\\
		0;|E|<2\end{cases}.\label{def:omega}\]
	
	We will let $\|M\|_{HS}=\Tr(MM^*)$ denote the Hilbert-Schmidt norm on $N$-by-$N$ matrices, and by $\|M\|=\sup_{x:\|x\|=1}\|Mx\|$, the operator norm of $M$. We define the $L^2$-Lipschitz norm on functions $f:\R\to \R$ as $\|f\|_{\mathrm{Lip}}=\sup_{x\neq y}\frac{|f(x)-f(y)|}{|x-y|}$. Given this, we lastly define the Wasserstein-1 distance on probability measures on $\R$ by
	\[\W_1(\mu,\nu)=\sup\{|\int_{\R}f(x)(\mu(dx)-\nu(dx))|:\|f\|_{\mathrm{Lip}}\le 1\}.\]
	
	\subsection{Acknowledgments} The author was partially supported by grants NSF DMS-1653552 and NSF DMS-1502632 while completing this project. The author would also like to thank their advisor, Antonio Auffinger, for introducing them to this problem and for many comments on earlier versions of this manuscript.
	\section{Outline of the Proof of Theorem \ref{Main-Theorem} and Auxiliary Results\label{section:outline}}
	
	In this section, we will introduce a series of lemmas (to be proved in later sections) from which we may derive Theorem \ref{Main-Theorem}. The starting point of our analysis will be an asymptotic upper bound for $\E[\Crit_{N,2}(B, I_O)]$.
	\begin{proposition}
		\label{prop:KR-2point}
		For $p,q\ge 5$, $0<\gamma<1$, nice $B\subseteq \R$, and nice $I_O\subseteq (-1,1)^2$, we have that
		\[\limsup_{N\to \infty}\frac{1}{N}\log\E[\Crit_{N,2}(B,I_O)]\le \sup_{E_1,E_2\in B,(r,t)\in I_O}\Sigma_{p,q,\gamma,2}(r,t,E_1,E_2).\]
	\end{proposition}
	
	As is standard in the calculation of the expectation of critical points \cite{pspin-one,pspin-second,bipartite}, this result will follow from careful analysis of the integrand appearing after an application of the Kac-Rice formula. The proof will be given in Section \ref{section:KR-2point}.
	
	We will also need the following result which verifies that the only overlaps in $\partial[-1,1]^2$ which contribute are in $\{\pm 1\}^2$.
	
	\begin{lemma}\label{lem:boundary}
		For $p,q\ge 5$, and $0<\gamma<1$, we have that \[\E[\Crit_{N,2}(\R,\partial [-1,1]^2-\{\pm 1\}^2)]=0.\]
	\end{lemma}
	
	We will now focus our analysis on the function $\Sigma_{p,q,\gamma,2}$. To understand this, we will first provide some results that will allow us to understand the supremum of $\Sigma_{p,q,\gamma,2}$ over $(r,t)$. Our first result will allow us to restrict our analysis to points where $E_1=E_2$.
	\begin{lemma}\label{lem:diagonal lemma}
		For $p,q\ge 5$, $0<\gamma<1$, nice $B\subseteq (-\infty,-E_{p,q,\gamma;\infty})$ and any $r,t\in (-1,1)$, we have that
		\[\sup_{E_1,E_2\in B}\Sigma_{p,q,\gamma,2}(r,t,E_1,E_2)\le \sup_{E\in B}\Sigma_{p,q,\gamma,2}(r,t,E,E).\]
		In particular we have that
		\[\limsup_{N\to\infty} \frac{1}{N}\log(\E[\Crit_{N,2}(B,(-1,1)^2)])\le \sup_{E\in B,r,t\in (-1,1)}\Sigma_{p,q,\gamma,2}(r,t,E,E).\]
	\end{lemma}
	
	The proof of this follows from similar arguments to the first part of Lemma 6 in \cite{pspin-second} and will be omitted.
	
	Let us define $\Psi_{p,q,\gamma}(r,t,E):=\Psi_{p,q,\gamma}(r,t,E,E)$ and similarly for $\Sigma_{2,p,q,\gamma}$. To complement this reduction, we will now also provide a result that identifies the supremum of $\Psi_{p,q,\gamma}$ for suitably small $E$. Let us denote $E_{p,q,\gamma;th}:=\sq{2(\gamma\log(p-1)+(1-\gamma)\log(q-1))}$.
	\begin{lemma}\label{lem:key-lemma}
		For $p,q\ge 10$, $0<\gamma<1$, and $|E|\le E_{p,q,\gamma;th}$ we have that
		\[\sup_{r,t\in (-1,1)}\Psi_{p,q,\gamma}(r,t,E)\le \Psi_{p,q,\gamma}(0,0,E)=2\Theta_{p,q,\gamma}(E).\label{eqn:key-lemma}\]
	\end{lemma}
	
	Now we reach the key technical results needed in the proof of Theorem \ref{Main-Theorem}. The primary estimate, which also is the limiting factor on improving the threshold of $p,q\ge 97$ is the following estimate which controls when $E_{p,q,\gamma;0}$ is smaller than the threshold required by Lemma \ref{lem:key-lemma}.
	\begin{lemma}\label{lem:E0-bound}
		For $p,q\ge 97$ and $0<\gamma<1$, we have that $E_{p,q,\gamma;0}<E_{p,q,\gamma;th}$.
	\end{lemma}
	
	Finally, we will require an additional estimate to address the case where $E>E_{p,q,\gamma;th}$. 
	
	\begin{lemma}\label{lem:negative bound}
		For $p,q\ge 97$, $0<\gamma<1$, and $r,t\in (-1,1)$ we have that $\Sigma_{2,p,q,\gamma}(E,r,t)\le 0$ for $|E|>E_{p,q,\gamma;th}$.
	\end{lemma}
	
	The proofs for these results will be given in Section \ref{section:lemma-section}. We mention that the condition $p,q\ge 97$ in Lemma \ref{lem:negative bound} is primarily needed as it takes Lemma \ref{lem:E0-bound} as input. Assuming that the bounds on Lemma \ref{lem:E0-bound} are improved, the proof of Lemma \ref{lem:negative bound} could likely be improved as well.
	
	With these preliminary results, we are ready to prove Theorem \ref{Main-Theorem}.
	
	\begin{proof}[Proof of Theorem \ref{Main-Theorem}]
		To begin, we will show that there is a choice of $\epsilon>0$ such that
		\[		\lim_{N\to\infty}\frac{1}{N}\log(\frac{\E[\Crit_N((-E_{p,q,\gamma;0}-\epsilon,E))^2]}{\E[\Crit_N((-E_{p,q,\gamma;0}-\epsilon,E))]^2})=0.\label{proof-eqn:main-theorem-1}
		\]
		Specifically, we will take $\epsilon=E_{p,q,\gamma;th}-E_{p,q,\gamma;0}$. We have that $\epsilon>0$ by Lemma \ref{lem:E0-bound}. We observe that by Jensen's inequality to show (\ref{proof-eqn:main-theorem-1}) it suffices to show that
		\[\limsup_{N\to\infty}\frac{1}{N}\log(\frac{\E[\Crit_N((-E_{p,q,\gamma;0}-\epsilon,E))^2]}{\E[\Crit_N((-E_{p,q,\gamma;0}-\epsilon,E))]^2})\le 0.\label{proof-eqn:main-theorem}\] 
		By Lemmas \ref{lem:diagonal lemma} and \ref{lem:key-lemma}, we see that for $E\in (-E_{p,q,\gamma;0}-\epsilon,-E_{p,q,\gamma;\infty})$
		\[\limsup_{N\to \infty}\frac{1}{N}\log(\E[\Crit_{N,2}((-E_{p,q,\gamma;0}-\epsilon,E),(-1,1)^2)])\le \sup_{t\in (-E_{p,q,\gamma;0}-\epsilon,E)}2\Sigma_{p,q,\gamma}(t).\]
		In light of Lemma \ref{lem:boundary}, we see that
		\[\E[\Crit_{N,2}((-E_{p,q,\gamma;0}-\epsilon,E),\partial[-1,1]^2)]=\E[\Crit_{N,2}((-E_{p,q,\gamma;0}-\epsilon,E),\{\pm 1\}^2)].\]
		We now observe that $\Crit_{N,2}(B,\{\pm 1\}^2)$ counts the number of critical points $(\pm \sigma,\pm \tau)$ (possibly with different signs), such that $\D H(\sigma,\tau)=0$ and $H(\pm \sigma,\pm\tau)\in NB$. In particular, $\Crit_{N,2}(B,\{\pm 1\}^2)\le 4\Crit_{N}(B)$. Combining these, we see that
		\[\limsup_{N\to \infty}\frac{1}{N}\log(\E[\Crit_{N,2}((-E_{p,q,\gamma;0}-\epsilon,E),[-1,1]^2)])\le \]\[\sup_{t\in (-E_{p,q,\gamma;0}-\epsilon,E)}[\max(\Sigma_{p,q,\gamma}(t),2\Sigma_{p,q,\gamma}(t))].\]
		Now we know by definition that $\Sigma_{p,q,\gamma}(t)$ is positive somewhere on $(-E_{p,q,\gamma;0},-E_{p,q,\gamma;0}+\delta)$ for each $\delta>0$. Thus we see that the supremum must be positive, so that
		\[\sup_{t\in (-E_{p,q,\gamma;0}-\epsilon,E)}[\max(\Sigma_{p,q,\gamma}(t),2\Sigma_{p,q,\gamma}(t))]=\sup_{t\in (-E_{p,q,\gamma;0}-\epsilon,E)}2\Sigma_{p,q,\gamma}(t),\]
		which in view of (\ref{Single point complexity}) completes the proof of (\ref{proof-eqn:main-theorem}).
		
		To complete the proof, it suffices to show in addition that
		\[\limsup_{N\to \infty}\frac{1}{N}\log(\E[\Crit_{N}((-\infty,-E_{p,q,\gamma;0}-\epsilon))^2])\le 0.\label{eqn:ignore-1249}\]
		Following the proof above, we see that the left-hand side of (\ref{eqn:ignore-1249}) is bounded above by
		\[\max(\sup_{u\le -E_{p,q,\gamma;0}-\epsilon,r,t\in (-1,1)} \Sigma_{2,p,q,\gamma}(r,t,u),\sup_{u\le -E_{p,q,\gamma;0}-\epsilon}\Sigma_{p,q,\gamma}(u)).\]
		The second term in this maximum is negative by the definition of $E_{p,q,\gamma;0}$, so we are left to show the same for the first term. Noting that $E_{p,q,\gamma;0}+\epsilon=E_{p,q,\gamma;th}$ though, this follows from Lemma \ref{lem:negative bound}.
	\end{proof}
	
	\section{Proof of Lemmas \ref{lem:boundary},  \ref{lem:key-lemma}, \ref{lem:E0-bound}, and \ref{lem:negative bound}\label{section:lemma-section}}
	In this section, we will prove some of the lemmas introduced in the last section. We will begin with the proof of Lemma \ref{lem:boundary}, for which we will first provide the following general result, which will essentially follow from Proposition 6.5 of \cite{azais}.
	\begin{lemma}\label{lem:azais-prelim}
		Let $T\subset M$ be a compact subset of a Riemannian manifold $M$, and $Z:M\to TM$ be a smooth Gaussian vector field on $M$. Let us assume that for each $t\in T$ we have that $Z(t)\in T_tM$ is a non-degenerate Gaussian vector. Then
		\[\P\left(\exists t\in T,\text{ such that }Z(t)=0\text{ and }\det(\D Z(t))=0\right)=0,\label{eqn:never-edit}\]
		where here $\D Z(t)$ denotes the Riemannian gradient of the vector field $Z$ evaluated at $t$.
	\end{lemma}
	\begin{proof}
		The first case of Proposition 6.5 of \cite{azais} is exactly this result when $M$ is the Euclidean space with the standard Euclidean metric. On the other hand, if $M$ is Euclidean with a nonstandard metric, then with respect to the standard basis Euclidean basis of $TM$, we have that
		\[[\D Z(t)]_{i,j}=\sum_{k}g^{j,k}(t)\partial_k Z_i(t),\]
		where $\partial_k$ denotes the $k$-th Euclidean partial derivative and $[g^{i,j}(t)]_{i,j}$ denotes the inverse metric tensor. In particular, we see that $\det(\D Z(t))=\det(g(t))^{-1}\det(\partial Z(t))$, where $\partial Z$ denotes the standard Euclidean gradient, so that the set (\ref{eqn:never-edit}) is independent of the choice of metric, which establishes the special case where $M$ is Euclidean.
		
		For the case of general $M$, note that we may take a finite number of Euclidean neighborhoods of $M$, say $(U_i)_{i=1}^m$, such that $T\subseteq \cup_{i=1}^m U_i$. Now, as $U_i$ is Euclidean, we may write it as an increasing union of compact subsets $(K_{i,j})_{j=1}^{\infty}$. By the union bound, we see that to show (\ref{eqn:never-edit}), it suffices to show for $i=1,\dots m$ and $j\in \N$ that
		\[\P(\exists t\in T\cap K_{i,j},\text{ such that }Z(t)=0\text{ and }\det(\D Z(t))=0)=0.\label{eqn:never-euc}\]
		As $T\cap K^{i,j}$ is compact in $U_i$ though, (\ref{eqn:never-euc}) follows from the Euclidean special case, so we are done.
	\end{proof}
	
	\begin{proof}[Proof of Lemma \ref{lem:boundary}]
		By employing that $H_{N}(-\sigma,\tau)=(-1)^pH_N(\sigma,\tau)$, as well as switching $p$ and $q$, we see that it suffices to show that $\E[\Crit_{N,2}(\R,\{1\} \times (-1,1))]=0$. By the monotone convergence theorem, it further suffices to show for each $\epsilon>0$ that $\E[\Crit_{N,2}(\R,\{1\} \times [-1+\epsilon,1-\epsilon])]=0$.
		If we define
		\[S_N(\epsilon)=\{(\sigma,\tau),(\sigma',\tau')\in S^{N_1-1}(\sq{N_1})\times S^{N_2-1}(\sq{N_2}):\frac{1}{N_2}|(\tau,\tau')|\le 1-\epsilon\},\]
		then we see that $\Crit_{N,2}(\R,\{1\} \times [-1+\epsilon,1-\epsilon])$ is exactly the number of zeros of the vector-field $Z(\sigma,\tau,\sigma',\tau')=(\D H_N(\sigma,\tau),\D H_N(\sigma,\tau'))$ that lie in $S_N(\epsilon)$. Note that $Z$ does not depend on $\sigma'$ at all, so that as $N_2>1$ we have that $\det(\D Z(t))=0$ for all $t\in S_N(\epsilon)$. In Remark \ref{remark:non-degen-new} (see also the definition (\ref{def-h-H})) we show that the law of $Z(t)$ is non-degenerate for $t\in S_N(\epsilon)$, so applying Lemma \ref{lem:azais-prelim} and the observation that the determinant vanishes everywhere, we see that
		\[\P(\Crit_{N,2}(\R,\{1\} \times [-1+\epsilon,1-\epsilon])>0)=0,\]
		which completes the proof.
	\end{proof}
	
	Before we begin the proof of Lemma \ref{lem:key-lemma}, we will prove a preliminary result that simplifies the function $b_{p,q}$, introduced in (\ref{better U-def}). For the convenience
	of the reader, we recall that
	\[b_{p,q}(r,t):=(1-r^{p}t^q)(1-r^{2p-2}t^{2q-2})+(p-1)(1-r^2)r^{p-2}t^{q}(1-r^{p}t^{q-2})\]\[+(q-1)(1-t^2)r^{p}t^{q-2}(1-r^{p-2}t^{q}).\]
	\begin{lemma}\label{key:b-simp}
		For $p,q\ge 2$ and $r,t\in [0,1)^2$, we have that
		\[b_{p,q}(r,t)\ge 2r^{p-2}t^{q-2}((p-1)(1-r^2)t^2+(q-1)(1-t^2)r^2)\frac{(1-r^{p-2}t^{q})(1-r^{p}t^{q-2})}{(1-r^{p-1}t^{q-1})}.\label{key:third-inequality}\]
	\end{lemma}
	\begin{proof}
		By continuity, it suffices to assume that $r,t\in (0,1)^2$. To prove (\ref{key:third-inequality}) we first show that it will suffice to prove the following statements 
		\[
		(1-r^{2p-2}t^{2q-2})\ge r^{p-2}t^{q-2}((p-1)(1-r^2)t^2+(q-1)(1-t^2)r^2)\label{key:first-inequality},
		\]
		\[
		(1-r^{p-2}t^q)+(1-r^{p}t^q)\ge 2\frac{(1-r^{p-2}t^{q})(1-r^{p}t^{q-2})}{(1-r^{p-1}t^{q-1})}.\label{key:second-inequality}
		\]
		Let us assume these statements for the moment. By applying (\ref{key:first-inequality}) to the first term of $b_{p,q}$ we get that
		\begin{gather}
			b_{p,q}(r,t)\ge (1-r^p t^q)r^{p-2}t^{q-2}((p-1)(1-r^2)t^2+(q-1)(1-t^2)r^2)+\\
			(p-1)(1-r^2)r^{p-2}t^{q}(1-r^{p}t^{q-2})+(q-1)(1-t^2)r^{p}t^{q-2}(1-r^{p-2}t^{q})=\\
			r^{p-2}t^{q-2}[(p-1)(1-r^2)t^2((1-r^{p}t^{q-2})+(1-r^{p}t^q))+\\
			(q-1)(1-t^2)r^2((1-r^{p-2}t^{q})+(1-r^{p}t^q))].\label{eqn:ignore-56}
		\end{gather}
		Applying (\ref{key:second-inequality}), and the same equality with $(p,r)$ and $(q,t)$ interchanged, we may further lower bound the right-hand side of (\ref{eqn:ignore-56}) by
		\[
		2r^{p-2}t^{q-2}((p-1)(1-r^2)t^2 +
		(q-1)(1-t^2)r^2 )\frac{(1-r^{p-2}t^{q})(1-r^{p}t^{q-2})}{(1-r^{p-1}t^{q-1})},
		\]
		which gives (\ref{key:third-inequality}).
		We now proceed with the proofs of (\ref{key:first-inequality}) and (\ref{key:second-inequality}), beginning with (\ref{key:first-inequality}). Let us denote
		\[j_{p,q}(r,t)=(1-r^{2p-2}t^{2q-2})-r^{p-2}t^{q-2}((p-1)(1-r^2)t^2+(q-1)(1-t^2)r^2),\]
		so that (\ref{key:first-inequality}) is equivalent to the $j_{p,q}(r,t)\ge 0$. We compute that
		\[
		r\frac{d}{dr}j_{p,q}(r,t)+t\frac{d}{dt}j_{p,q}(r,t)=
		-2(p+q-2)r^{2p-2}t^{2q-2}-(p+q-2)(p-1)r^{p-2}t^q+\]
		\[
		(p+q)(p-1)r^p t^q-(p+q-2)(q-1)r^{p}t^{q-2}+(p+q)(q-1)r^p t^q=
		\]
		\[
		(p+q-2)r^{p-2}t^{q-2}\bar{j}_{p,q}(r,t),\]
		\[
		\bar{j}_{p,q}(r,t):=-2r^p t^q+(p+q)r^2t^2-(p-1)t^2-(q-1)r^2.\label{eqn:derivative j}
		\]
		We will first prove that $\bar{j}_{p,q}(r,t)\le 0$. To do this, note that
		\[\frac{d}{dp}\bar{j}_{p,q}(r,t)=-2r^p t^q\log(r)-(1-r^2)t^2\le -2r^2t^2\log(r)-(1-r^2)t^2.\]
		We have that
		\[\frac{d}{dr}(-2r^2t^2\log(r)-(1-r^2)t^2)=-4r t^2 \log(r)>0,\]
		and as this function vanishes when $r=1$, we see that
		\[\frac{d}{dp}\bar{j}_{p,q}(r,t)\le -2r^2t^2\log(r)-(1-r^2)t^2\le 0.\]
		By symmetry we also see that $\frac{d}{dq}\bar{j}_{p,q}(r,t)\le 0$. In particular 
		\[\bar{j}_{p,q}(r,t)\le \bar{j}_{2,2}(r,t)=-r^2-t^2+2r^2t^2\le 0.\]
		Now that we have shown that $\bar{j}_{p,q}(r,t)\le 0$ we see that for $r,t\in (0,1)$
		\[r\frac{d}{dr}j_{p,q}(r,t)+t\frac{d}{dt}j_{p,q}(r,t)\le 0.\] 
		Thus in particular for $\lam\in (0,\min(r^{-1},t^{-1}))$ we see that $\frac{d}{d\lam}j_{p,q}(\lam r, \lam t)\le 0$, so taking $\lam=\min(r^{-1},t^{-1})$ we see that it suffices to show that $j_{p,q}(1,t)\ge 0$ and $j_{p,q}(r,1)\ge 0$. As $j_{p,q}(r,t)=j_{q,p}(t,r)$, it additionally suffices to only show the first statement. In this case
		\[j_{p,q}(1,t)=(1-t^{2q-2})-(q-1)t^{q-2}(1-t^2).\]
		Now we note that by the AM-GM inequality, we have that
		\[\frac{1}{q-1}\frac{(1-t^{2q-2})}{1-t^2}=\frac{1}{q-1}\sum_{i=0}^{q-2}t^{2i}\ge t^{q-2}.\]
		Rearranging this exactly gives that $j_{p,q}(1,t)\ge 0$, and completes the proof of (\ref{key:first-inequality}).
		
		We now will show (\ref{key:second-inequality}). We observe that
		\[(1-r^{p-1}t^{q-1})((1-r^{p-2}t^q)+(1-r^{p}t^q))-2(1-r^{p-2}t^{q})(1-r^{p}t^{q-2})=\]\[r^{p-2}t^{q-2}\bar{l}_{p,q}(r,t),\;\;\bar{l}_{p,q}(r,t):=((t-r)^2+r^2(1-t^2)+r^{p-1}t^{q+1}+r^{p+1}t^{q+1}-2r^p t^q).\label{eqn:second in}\]
		Thus by clearing the denominator, we see that (\ref{key:second-inequality}) is equivalent to showing that $\bar{l}_{p,q}(r,t)\ge 0$.
		
		We note that $\lim_{p\to \infty}\bar{l}_{p,q}(r,t)=(t-r)^2+r^2(1-t^2)$. In addition, as 
		\[\frac{d}{dp}\bar{l}_{p,q}(r,t)=r^{p-1} t^q(t+r^2t-2r)\log(r),\]
		we see that at a fixed point $(r,t)$, the sign of $\frac{d}{dp}\bar{l}_{p,q}(r,t)$ does not depend on $p$. In particular, we see that \[\bar{l}_{p,q}(r,t)\ge  \min(\bar{l}_{1,q}(r,t),(t-r)^2+r^2(1-t^2)).\]
		The same argument in $q$ implies that 
		\[\bar{l}_{1,q}(r,t)\ge  \min(\bar{l}_{1,1}(r,t),(t-r)^2+r^2(1-t^2)).\]
		Lastly noting that $\bar{l}_{1,1}(r,t)=2(r-t)^2$, we see that both of these are clearly positive, so we see that $\bar{l}_{p,q}(r,t)\ge 0$, which establishes (\ref{key:second-inequality}).
	\end{proof}
	
	\begin{proof}[Proof of Lemma \ref{lem:key-lemma}]
		Let us denote \[k_{p,q}(r,t)=\frac{t^qr^p(1-r^{p-2}t^q)(1-r^p t^{q-2})}{b_{p,q}(r,t)},\label{eqn:k-def}\]
		\[Q_{p,q,\gamma}(r,t,E)=g_{p,q,\gamma}(r,t)+E^2k_{p,q}(r,t).\]
		In sight of (\ref{eqn:psi-def}) and (\ref{better U-def}), we see that $Q_{p,q,\gamma}$ consists of all the terms in $\Sigma_{p,q,\gamma,2}$ which depend on $(r,t)$, and so that to establish (\ref{eqn:key-lemma}), it suffices to show that for $r,t\in (-1,1)$ and $|E|\le E_{p,q,\gamma;th}$
		\[Q_{p,q,\gamma}(r,t,E)\le 0.\]
		We will now show that it suffices to assume that $r,t\in (0,1)$. By continuity, we only need to show that for $r,t\in (-1,1)-\{0\}$ we have that
		\[Q_{p,q,\gamma}(r,t,E)\le Q_{p,q,\gamma}(|r|,|t|,E).\label{eqn:ignore-scene}\]
		
		To begin we observe that (\ref{eqn:b-det-Y}) below expresses that $p^{-1}q^{-1}b_{p,q}(r,t)$ is the determinant of covariance matrix of a non-degenerate Gaussian vector (\ref{eqn:ignore-bing}), so that we must have that $b_{p,q}(r,t)>0$ for all $r,t\in (-1,1)$. Now note that $g_{p,q,\gamma}(r,t)$ is even in both $r$ and $t$. We also observe that if $\sign(r)^p\sign(t)^q=1$ then $k_{p,q}(r,t)=k_{p,q}(|r|,|t|)$, so that $Q_{p,q,\gamma}(r,t,E)=Q_{p,q,\gamma}(|r|,|t|,E)$. If $\sign(r)^p\sign(t)^q=-1$, we see that $k_{p,q}(r,t)$ is negative, so that
		\[Q_{p,q,\gamma}(r,t,E)\le g_{p,q,\gamma}(r,t)\le Q_{p,q,\gamma}(|r|,|t|,E).\]
		Together these verify (\ref{eqn:ignore-scene}).
		
		We will now assume that $r,t\in (0,1)$ for the remainder of the proof. By Lemma \ref{key:b-simp}, we see that if we denote
		\[m_{p,q}(r,t)=(p-1)(1-r^2)t^2+(q-1)(1-t^2)r^2,\]
		\[\bar{k}_{p,q}(r,t)=\frac{1}{2}\frac{t^2r^2(1-r^{p-1}t^{q-1})}{m_{p,q}(r,t)},\label{eqn:ignore-622}\]
		then we have that $\bar{k}_{p,q}(r,t)\ge k_{p,q}(r,t)$. Let us define
		\[\bar{Q}_{p,q,\gamma}(r,t,E):=g_{p,q,\gamma}(r,t)+E^2\bar{k}_{p,q}(r,t).\]
		As we have that $\bar{Q}_{p,q,\gamma}(r,t,E)\ge Q_{p,q,\gamma}(r,t,E)$ it will suffice to verify the corresponding statements for $\bar{Q}_{p,q,\gamma}$. Now as $\bar{k}_{p,q}$ is positive we see that for $|E|\le |E'|$
		\[\bar{Q}_{p,q,\gamma}(r,t,E)\le \bar{Q}_{p,q,\gamma}(r,t,E').\]
		As we have in addition that $\bar{Q}_{p,q,\gamma}(r,t,-E)=\bar{Q}_{p,q,\gamma}(r,t,E)$, we see that to prove the desired statement it is sufficient to show that $\bar{Q}_{p,q,\gamma}(r,t,E_{p,q,\gamma;th})\le 0$. On the other hand, we observe that
		\[\bar{Q}_{p,q,\gamma}(r,t,E_{p,q,\gamma;th})=\gamma \bar{Q}_{p,q,1}(r,t,E_{p,q,1;th})+(1-\gamma)\bar{Q}_{p,q,0}(r,t,E_{p,q,0;th}).\]
		As $\bar{Q}_{p,q,0}(r,t,E_{p,q,0;th})=\bar{Q}_{q,p,1}(t,r,E_{q,p,1;th})$, we see that it is also sufficient to assume that $\gamma=1$. Let us denote $\bar{Q}_{p,q}(r,t):=\bar{Q}_{p,q,1}(r,t,E_{p,q,1;th})$.
		
		We will first show that $\frac{d}{dt}\bar{k}_{p,q}(r,t)\ge 0$. Direct calculation yields that
		\[\frac{d}{dt}\bar{k}_{p,q}(r,t)=\frac{(q-1)r^{p+1}t^q}{2m_{p,q}(r,t)^2}\hat{j}_{p,q}(r,t)\]\[\hat{j}_{p,q}(r,t)=2r^{3-p}t^{1-q}-2r^2-m_{p,q}(r,t).\]
		Thus to show that $\frac{d}{dt}\bar{k}_{p,q}(r,t)\ge 0$, it suffices to show that $\hat{j}_{p,q}(r,t)\ge 0$. To do this we first show that for $p\ge 3$ and $q\ge 1$, that we have that $\frac{d}{dq}\hat{j}_{p,q}(r,t)\ge 0$. We compute that
		\[\frac{d}{dq}\hat{j}_{p,q}(r,t)=-r^2+r^2t^2-2r^{3-p}t^{1-q}\log(t);\;\;\frac{d^2}{dq^2}\hat{j}_{p,q}(r,t)=2r^{3-p}t^{1-q}\log(t)^2,\]
		\[\frac{d^2}{dqdp}\hat{j}_{p,q}(r,t)=2r^{3-p}t^{1-q}\log(r)\log(t).\]
		In particular, we see that $\frac{d^2}{dq^2}\hat{j}_{p,q}(r,t)\ge 0$ and $\frac{d^2}{dqdp}\hat{j}_{p,q}(r,t)\ge 0$, so that \[\frac{d}{dq}\hat{j}_{p,q}(r,t)\ge \frac{d}{dq}\hat{j}_{p,q}(r,t)\big|_{(p,q)=(3,1)}=-r^2(1-t^2)-2\log(t)\ge -(1-t^2)-2\log(t)\ge 0,\]
		as the final expression may easily be checked to be decreasing on $(0,1)$ and vanishes at $1$. We note that $\hat{j}_{p,q}(r,t)=\hat{j}_{q+2,p-2}(t,r)$ so that we also have that $\frac{d}{dp}\hat{j}_{p,q}(r,t)\ge 0$ for $p\ge 3$ and $q\ge 1$. Combining these, we see that
		\[\hat{j}_{p,q}(r,t)\ge \hat{j}_{3,1}(r,t)=2(1-r^2)(1-t^2)\ge 0.\]
		This completes the proof that $\frac{d}{dt}\bar{k}_{p,q}(r,t)>0$. By symmetry, we see that $\frac{d}{dr}\bar{k}_{p,q}(r,t)>0$.
		\noindent
		Now we observe that
		\[\frac{d}{dt}g_{p,q,1}(r,t)=\frac{qr^{2p-2}t^{2q-1}}{1-r^{2p-2}t^{2q}}>0.\]
		Combining these observations, we see that
		\[\frac{d}{dt}\bar{Q}_{p,q}(r,t)>0.\]
		It is clear that $\bar{Q}_{p,q}(r,t)$ extends continuously to $(r,t)\in [0,1)\times [0,1]$, and so we see that
		\[\sup_{r,t\in (0,1)^2}\bar{Q}_{p,q}(r,t)\le \sup_{r\in (0,1)}\bar{Q}_{p,q}(r,1).\]
		Thus we are reduced to showing that for $r\in (0,1)$,
		\[\bar{Q}_{p,q}(r,1)=\frac{1}{2}\log(\frac{1-r^2}{1-r^{2p-2}})+\frac{\log(p-1)}{p-1}\frac{1-r^{p-1}}{1-r^2}r^2\le 0.\]
		We will denote $\bar{Q}_{p}(r):=\bar{Q}_{p,q}(r,1)$. Observing that $\frac{1-r^{2p-2}}{1-r^2}=\sum_{i=0}^{p-2}r^{2i}\ge 1+r^2$, we see that
		\[\bar{Q}_{p}(r)\le -\frac{1}{2}\log(1+r^2)+\frac{\log(p-1)}{p-1}\frac{1}{1-r^2}r^2=:\hat{Q}_{p}(r).\]
		As $\log(p-1)/(p-1)$ is decreasing in $p>4$, we see that $\frac{d}{dp}\hat{Q}_p(r)<0$ for $p>4$. We observe that for $r\in [0,0.61)$ \[\log(1+r^2)\ge \frac{\log(1+0.61^2)}{0.61^2}r^2.\] From this we see that for $r\in [0,0.61)$
		\[\hat{Q}_{10}(r)\le (-\frac{\log(1+0.61^2)}{0.61^2}+\frac{\log(10-1)}{10-1}\frac{1}{1-0.61^2})r^2=-(0.461\dots)r^2\le 0.\]
		Thus for $p\ge 10$ and $r\in [0,0.61)$, we have that $\bar{Q}_{p}(r)\le 0$. Now we recall that for $r\in (0,1)$, the function $(1-r^a)/a$ is decreasing in $a>0$. In particular, we see that
		\[\bar{Q}_{p}(r)\le \frac{1}{2}\log(\frac{1-r^2}{1-r^{2p-2}})+\frac{\log(p-1)}{p-2}\frac{1-r^{p-2}}{1-r^2}r^2.\]
		The right-hand side of this inequality coincides with the function $\tilde{Q}_p(r)$ introduced on pg. 3415 of \cite{pspin-second}. On pg. 3416 it is moreover shown that for $p\ge 10$ and $r\in [0.6,1)$, we have that $\tilde{Q}_p(r)<0$. Together these show that $\bar{Q}_p(r)\le 0$ for $p\ge 10$, completing the proof.
	\end{proof}
	
	\begin{remark}
		One may check numerically that $\bar{Q}_{p}(r)\le 0$ for $r\in (0,1)$ and $p\ge 3$. Given this claim, the proof of Lemma \ref{lem:key-lemma} in fact shows that Lemma \ref{lem:key-lemma} still holds under the weaker assumption that $p,q\ge 3$. As this improvement does not affect our main result, though, we choose not to pursue it.
	\end{remark}

	We will now begin to address the proof of Lemma \ref{lem:E0-bound}. As a preliminary step, we will obtain an upper bound on $\Sigma_{p,q,\gamma}(E)$. This bound is essentially equivalent to the one obtained in Section 4.1 of \cite{tuca} for a single point. We will say a symmetric random matrix $M_n$ of size $n$-by-$n$ is a $\GOE(n)$ matrix if its entries above and on the main diagonal are independent centered Gaussian random variables with variance $\E[(M_n)_{ij}^2]=\frac{1}{n}(1+\delta_{ij})$. It is well-known (see, for example \cite{pspin-one} or \cite{exponential}) that
	\[\lim_{N\to\infty}\frac{1}{N}\log(\E[|\det(M_{N}-EI)|])=\Omega_{sc}(E),\]
	uniformly in compact subsets of $E\in \R$. Now we recall Fischer's Inequality: For $M$ positive semi-definite, we have that
	\[M=\begin{bmatrix}A&B\\B^t&C \end{bmatrix}\implies \det(M)\le \det(A)\det(C).\label{fischers-inequality}\]
	We will apply this to the formula for $\Sigma_{p,q,\gamma}$ derived in Lemma \ref{lem:mckenna lem:G-det upperbound}. Let us denote $\gamma_{1}:=p^{-1}(p-1)(N_1-1)/N$ and $\gamma_{2}:=q^{-1}(q-1)(N_2-1)/N$. Applying (\ref{fischers-inequality}) we see that
	\[\E[|\det(H_N^D-EI)|I(H_N^D\ge EI)]\le\]\[ \prod_{i=1}^2\gamma_i^{(N_i-1)/2}\E[|\det(M_{N_i-1}-\gamma_i^{-1/2}EI)|],\label{eqn:ignore-1129}\]
	where $M_n$ denotes a $\GOE(n)$ matrix as above. In particular, we see that for $E<-E_{p,q,\gamma;\infty}$
	\[\Sigma_{p,q,\gamma}(E)\le C_{p,q,\gamma}+
	\gamma\Omega_{sc}(E\sq{\frac{p}{(p-1)\gamma}})+(1-\gamma)\Omega_{sc}(E\sq{\frac{q}{(q-1)(1-\gamma)}})-\frac{1}{2}E^2.\label{tuca-upper bound}\]
	This upper bound will serve as our method to bound $E_{p,q,\gamma;0}$. 
	
	\begin{proof}[Proof of Lemma \ref{lem:E0-bound}]
		It suffices to show that for $E\le -E_{p,q,\gamma;th}$, we have that $\Sigma_{p,q,\gamma}(E)<0$. To do this we first observe that by (\ref{tuca-upper bound}), it suffices to show that the right-hand side of (\ref{tuca-upper bound}) is negative when $E\le -E_{p,q,\gamma;th}$.
		
		To simplify this expression, first note that by the inequality \[\gamma\log(\gamma)+(1-\gamma)\log(1-\gamma)\ge -\log(2),\]
		we have that 
		\[C_{p,q,\gamma}\le \frac{1}{2}[1+\log(2)+\gamma \log(p)+(1-\gamma)\log(q)].\]
		Let us denote $c_1=\sq{96/97}=0.994...$.
		Noting that $\Omega_{sc}(E)$ is monotonically decreasing for $E\in (-\infty,0)$ we see that the right-hand side of (\ref{tuca-upper bound})  is upper bounded for $E<0$ by
		\[\bar{\Sigma}_{p,q,\gamma}(E):=\frac{1}{2}[1+\log(2)+\gamma \log(p)+(1-\gamma)\log(q)]+\]\[\gamma \Omega_{sc}(Ec_1^{-1}/\sq{\gamma})+(1-\gamma) \Omega_{sc}(Ec_1^{-1}/\sq{1-\gamma})-\frac{1}{2}E^2.\label{eqn:ignore-645}\]
		This expression is even in $E$, so it suffices to show that for $E>E_{p,q,\gamma;th}$, we have that $\bar{\Sigma}_{p,q,\gamma}(E)<0$. We define the function
		\[\hat{\Sigma}_{p,q,\gamma}(E)=\bar{\Sigma}_{p,q,\gamma}(E)+\frac{E^2}{4}-\frac{E_{p,q,\gamma;th}^2}{4}.\]
		We note that for $E\ge E_{p,q,\gamma;th}$ we have that $\bar{\Sigma}_{p,q,\gamma}(E)\le \hat{\Sigma}_{p,q,\gamma}(E)$ so it suffices to show that $\hat{\Sigma}_{p,q,\gamma}(E)<0$. To do so we will first show that $\hat{\Sigma}_{p,q,\gamma}'(E)\le 0$ when $E>3c_1=2.984...$. 
		
		We first note that for $E>2$
		\[\Omega'_{sc}(E)=\int \frac{1}{E-x}\frac{\sq{4-x^2}}{2\pi}dx\le \frac{1}{E-2}.\]
		If we denote
		\[f(\gamma,E)=\frac{\gamma}{E/\sq{\gamma}-2}+\frac{(1-\gamma)}{E/\sq{1-\gamma}-2},\]
		we see from this bound that
		\[\hat{\Sigma}_{p,q,\gamma}'(E)\le f(\gamma,c_1^{-1}E)-E/2.\]
		We that $f(\gamma,c_1^{-1}E)-E/2$ is decreasing in $E$ for $E\ge 2c_1$, so that it suffices to show that $f(\gamma,3)<3c_1/2$. On the other hand, one may verify by direct computation that $f(\gamma,3)$ is convex in $\gamma\in (0,1)$, so that $f(\gamma,3)\le \max(f(0,3),f(1,3))=1<3c_1 /2$. This establishes the desired claim for the derivative.
		
		Now noting that $E_{p,q,\gamma;th}\ge E_{97,97,\gamma;th}=3.021...>3c_1$, we only need to show that $\hat{\Sigma}_{p,q,\gamma}(E_{97,97,\gamma;th})<0$. 		
		We note that 
		\[\gamma\log(p/(p-1))+(1-\gamma)\log(q/(q-1))\le \log(97/96),\]
		so that we see that
		\[\hat{\Sigma}_{p,q,\gamma}(E_{97,97,\gamma;th})\le \gamma \Omega_{sc}(Ec_1^{-1}/\sq{\gamma})+(1-\gamma) \Omega_{sc}(Ec_1^{-1}/\sq{(1-\gamma)})\]\[-\frac{E_{97,97,\gamma;th}^2}{4}+\frac{1}{2}(1+\log(2)+\log(97/96)).\label{eqn:ignore-617}\]
		To analyze the terms on the right-hand side, we note that for $E>2$
		\[\Omega_{sc}'(E)-E\Omega_{sc}''(E)=\frac{2}{\sq{E^2-4}}>0.\]
		In particular, we see that for $x\in (0,1)$,  
		\[\frac{d^2}{dx^2}(x\Omega_{sc}(E/\sq{x}))=-\frac{E}{4x^{3/2}}(\Omega_{sc}'(E/\sq{x})-(E/\sq{x})\Omega_{sc}''(E/\sq{x}))<0.\]
		This shows that $x\Omega_{sc}(E/\sq{x})$ is concave for $x\in (0,1)$. From this, we see that the function
		\[\bar{f}(\gamma):=\gamma \Omega_{sc}(Ec_1^{-1}/\sq{\gamma})+(1-\gamma) \Omega_{sc}(Ec_1^{-1}/\sq{(1-\gamma)})\]
		is concave on $(0,1)$, but as $\bar{f}(\gamma)=\bar{f}(1-\gamma)$, we see that $\sup_{\gamma\in (0,1)}\bar{f}(\gamma)\le \bar{f}(1/2)$. Thus expanding upon (\ref{eqn:ignore-617}) we see that $\hat{\Sigma}_{p,q,\gamma}(E_{97,97,\gamma;th})$ is less than
		\[\Omega_{sc}(E_{97,97,\gamma;th}c_1^{-1}/\sq{2})+\frac{E_{97,97,\gamma;th}^2}{4}+\frac{1}{2}(1+\log(2)+\log(97/96))=-0.0016...\]
	\end{proof}
	
	The remainder of this section will consist of the proof of Lemma \ref{lem:negative bound}. This proof will involve the same basic tools as the proof Lemma \ref{lem:E0-bound}, and in particular, crucially uses (\ref{tuca-upper bound}) as well.
	
	\begin{proof}[Proof of Lemma \ref{lem:negative bound}]
		As the function is even we may assume that $E$ is positive. Employing (\ref{eqn:k-def}), (\ref{tuca-upper bound}), and (\ref{better U-def}) we see that for $E>E_{p,q,\gamma;\infty}$
		\[\Sigma_{2,p,q,\gamma}(E,r,t)=2\Sigma_{p,q,\gamma}(E)-E^2+E^2k_{p,q}(r,t)\le \]
		\[2C_{p,q,\gamma}+2\gamma\Omega_{sc}(E\sq{\frac{p}{(p-1)\gamma}})+2(1-\gamma)\Omega_{sc}(E\sq{\frac{q}{(q-1)(1-\gamma)}})-E^2+E^2k_{p,q}(r,t).\]
		Let us denote the function on the right-hand side as $\hat{\Sigma}_{2,p,q,\gamma}(E,r,t)$. As this function differs from $\Sigma_{2,p,q,\gamma}$ by a term depending only on $E$, Lemma \ref{lem:key-lemma} shows that for $E=E_{p,q,\gamma;th}$ we have that
		\[\sup_{r,t\in (-1,1)}\hat{\Sigma}_{2,p,q,\gamma}(E,r,t)=\]\[2C_{p,q,\gamma}+2\gamma\Omega_{sc}(E_{p,q,\gamma;th}\sq{\frac{p}{(p-1)\gamma}})+2(1-\gamma)\Omega_{sc}(E_{p,q,\gamma;th}\sq{\frac{q}{(q-1)(1-\gamma)}}).\]
		The proof of Lemma \ref{lem:E0-bound} shows that the quantity on the right-hand side is negative, so in particular it suffices to show that for all $r,t\in (-1,1)$ the function $\hat{\Sigma}_{p,q,\gamma}(E,r,t)$ is decreasing in $E>E_{p,q,\gamma;th}$.
		
		In the course of the proof of Lemma \ref{lem:key-lemma}, we showed that $k_{p,q}(r,t)\le k_{p,q}(|r|,|t|)$ and in (\ref{eqn:ignore-622}) we defined a function $\bar{k}_{p,q,\gamma}(r,t)$ such that for $p,q\ge 3$, and $r,t\in (0,1)$, $k_{p,q,\gamma}(r,t)\le \bar{k}_{p,q,\gamma}(r,t)$. We showed in addition that $\frac{d}{dt}\bar{k}_{p,q,\gamma}(r,t),\frac{d}{dr}\bar{k}_{p,q,\gamma}(r,t)\ge 0$. In particular, for $r,t\in (-1,1)$ we have that $k_{p,q}(r,t)\le \lim_{\lam\to 1} \bar{k}_{p,q,\gamma}(\lam,\lam)=1/4$.
		
		Note as well that for $p>3$, one may check that $E_{p,q,\gamma;th}>\sq{2\gamma\log(p-1)}>2\sq{(p-1)\gamma/p}$. In particular, it suffices to show that for all $0<\gamma<1$ the that function
		\[\hat{h}_\gamma(E):=2\gamma\Omega_{sc}(E\sq{\frac{p}{(p-1)\gamma}})-\frac{3}{8}E^2\]
		is decreasing in $E$ for $E>2\sq{(p-1)\gamma/p}$. We note that $\Omega_{sc}(x)$ is strictly convex in $x$ for $x>2$ and that $\Omega_{sc}'(2)=1$. In particular, we have that
		\[\hat{h}_\gamma'(E)<\hat{h}_\gamma'\left(2\sq{\frac{(p-1)\gamma}{p}}\right)=\sq{\gamma}\left(2\sq{\frac{p}{(p-1)}}-3\sq{\frac{(p-1)}{p}}\right).\label{eqn:ignore-191}\]
		The function $2\alpha^{1/2}-3\alpha^{-1/2}$ easily checked to be negative for $\alpha<3/2$ so that the right-hand side of (\ref{eqn:ignore-191}) is negative, which verifies the claim.
		
	\end{proof}

	\section{Proof of Proposition \ref{prop:KR-2point}\label{section:KR-2point}}
	
	This section is devoted to the proof of Proposition \ref{prop:KR-2point}. First we will perform a normalization on $H_N$ to simplify the factors which appear in our analysis. For each $\ell$, let $S^{\ell}=S^{\ell}(1)=\{x\in \R^\ell:\|x\|=1\}$. Furthermore, let us define for $(\sigma,\tau)\in S^{N_1-1}\times S^{N_2-1}$  \[h_N(\sigma,\tau)=h_{p,q,\gamma,N}(\sigma,\tau)=\frac{1}{\sq{N}}H_{p,q,\gamma,N}(\sq{N_1}\sigma,\sq{N_2}\tau).\label{def-h-H}\]
	The function $h_{N}$ is a centered Gaussian random field with covariance
	\[\E[h_N(\sigma,\tau)h_N(\sigma',\tau')]=(\sigma,\sigma')^p(\tau,\tau')^q.\]
	The proof of Proposition \ref{prop:KR-2point} begins, as in \cite{pspin-one},\cite{pspin-second} and \cite{bipartite}, by an application of the Kac-Rice formula. To state this, we will need some prerequisite notation. Let us endow the product of spheres $S^{N_1-1}\times S^{N_2-1}$ with the standard Riemannian metric inherited from $\R^N$. Given a choice of a (piecewise) smooth orthonormal frame field, $(E_i)_{i=1}^{N-2}$ on $S^{N_1-1}\times S^{N_2-1}$, we define
	\[\D h_N(\sigma)=((E_i h_N)(\sigma))_{i=1}^{N-2},\;\;\;\; \D^2 h_N(\sigma)=((E_i E_j h_N)(\sigma))_{i,j=1}^{N-2}.\]
	\begin{lemma}\label{lem:KR-basic}
		For arbitrary $(E_i)_{i=1}^{N-2}$ as above, and any choice of nice $B\subseteq \R$ and nice $I_O\subseteq (-1,1)^2$, we have that
		\begin{gather}
			\E[\Crit_{N,2}(B,I_O)]=C_{N}\int_{I_O}(1-r^2)^{(N_1-3)/2}(1-t^2)^{(N_2-3)/2}\times\\
			\varphi_{\D h_N(\n),\D h_N(\n(r,t))}(0,0)\E[|\det(\D^2h_N(\n))\det(\D^2h_N(\n(r,t)))|\times\\
			I(h_N(\n),h_N(\n(r,t))\in \sq{N}B)|\D h_N(\n)=\D h_N(\n(r,t))=0]drdt
		\end{gather}
		where here \[C_N=C_{N_1,N_2}=\vol(S^{N_1-1})\vol(S^{N_2-1})\vol(S^{N_1-2})\vol(S^{N_2-2}),\] the term $\varphi_{\D h_N(\n),\D h_N(\n(r,t))}(0,0)$ denotes the density of the random vector $$(\D h_N(\n),\D h_N(\n(r,t)))$$ at the event $(\D h_N(\n)=\D h_N(\n(r,t))=0)$, and $\n(r,t)=(\n_{N_1}(r), \n_{N_2}(t))$ where \[\n_{N}(s)=(s,\sq{1-s^2},0,\dots, 0)\in S^{N-1}\]
		and $\n=\n(1,1)$.
	\end{lemma}
	
	The proof of this result will be given in Appendix \ref{section:appendix-KR}. To make further sense of this quantity, we will need to recall the structure of the correlations between $h_N$ and its derivatives.
	
	\begin{lemma}\label{lem:Gradient-covariance-statement}
		For $r,t\in (-1,1)$, there exists a choice of $(E_i)_{i=1}^{N-1}$ such that the density of $(\D h_N(\n),\D h_N(\n(r,t)))$ at the event $(\D h_N(\n)=0,\D h_N(\n(r,t))=0)$ is given by
		\[\varphi_{\D h_N(\n),\D h_N(\n(r,t))}(0,0)=\]
		\[(2\pi)^{-(N-2)}(p^2(1-t^{2q}r^{2p-2}))^{-(N_1-2)/2}(q^2(1-r^{2p}t^{2q-2}))^{-(N_2-2)/2}f_L(r,t),\label{eqn:graident-to-find}\]
		where $f_L(r,t)$ is a fixed continuous function, positive on $(-1,1)^2$, given by (\ref{appendix f_L}). In addition, the law of $(h_N(\n),h_N(\n(r,t)))$, conditioned on $(\D h_N(\n)=0,\D h_N(\n(r,t))=0)$, is a centered Gaussian vector with covariance matrix $\Sigma_U(r,t)$ given by (\ref{appendix Sigma_E}).
	\end{lemma}
	
	The proof of this result will be given in Appendix \ref{section:covariance-section}. Before starting our next result, we let $D_{N}$ denote the $(N-2)$-by-$(N-2)$ matrix with $(N_1-1)$ and $(N_2-1)$ blocks given by
	\[D_{N}=N^{1/4}\begin{bmatrix}\sq{p}I_{N_1-1}&0\\ 0&\sq{q}I_{N_2-1} \end{bmatrix}.\]
	In addition, let us denote by $e_{ij}^{n,m}$ the elementary $n$-by-$m$ matrix with entries $[e_{ij}^{n,m}]_{kl}=\delta_{ik}\delta_{jl}$, and denote $e_{ij}=e^{n,m}_{ij}$ when the choice of $n$ and $m$ is clear from context.
	\begin{lemma}\label{lem:Hessian-covariance-statement}
		For $r,t\in (-1,1)$, and with the same choice of $(E_i)_{i=1}^{N-1}$ as in Lemma \ref{lem:Gradient-covariance-statement}, the following holds. Conditional on \[(\D h_N(\n)=\D h_N(\n(r,t))=0,h_N(\n)=\sq{N}E_1,h_N(\n(r,t))=\sq{N}E_2),\] let us denote the law of \[(D^{-1}_{N}\D^2h_N(\n)D^{-1}_{N},D^{-1}_{N}\D^2h_N(\n(r,t))D^{-1}_{N})\]
		as  \[(M^{1}_{N}(r,t,E_1,E_2),M^{2}_{N}(r,t,E_1,E_2)).\] Then for $k=1,2$, we may write $M^{k}_{N}(r,t,E_1,E_2)$ in terms of square $(N_1-1)$ and $(N_2-1)$ blocks as
		\[M^{k}_{N}(r,t,E_1,E_2)=\begin{bmatrix}
			\hat{M}^k_{N,1,1}(r,t)&\hat{M}^k_{N,1,2}(r,t)\\
			\hat{M}^k_{N,2,1}(r,t)&\hat{M}^k_{N,2,2}(r,t)\end{bmatrix}\]\[-E_kI+\begin{bmatrix}
			m_{11}^k(r,t,E_1,E_2)e_{1,1}&m_{12}^k(r,t,E_1,E_2)e_{1,1}\\
			m_{21}^k(r,t,E_1,E_2)e_{1,1}&m_{22}^k(r,t,E_1,E_2)e_{1,1}\end{bmatrix},\]
		where $m_{ij}^k$ are given by fixed (deterministic) functions only dependant on $p,q$ and such that $m_{ij}^k(r,t,E_1,E_2)=m_{ji}^k(r,t,E_1,E_2)$ and $\hat{M}^k_{N,i,j}$ have centered Gaussian entries and satisfy $\hat{M}^k_{N,i,j}(r,t)=(\hat{M}^k_{N,j,i}(r,t))^t$. 
		
		Let us further decompose $\hat{M}^k_{N,i,j}(r,t)$ as a block matrix
		\[\hat{M}^k_{N,i,j}(r,t)=\begin{bmatrix}
			Q^k_{N,ij}(r,t)& V^k_{N,ji}(r,t)^t\\
			V^k_{N,ij}(r,t)&G^k_{N,ij}(r,t)\\
		\end{bmatrix},\]
		with dimensions so that $Q^k_{N,ij}(r,t)$ is of size $1$-by-$1$ and $G^{k}_{N,ij}(r,t)$ is of size $(N_j-2)$-by-$(N_i-2)$. Then for $k=1,2$, the entries satisfying the following:
		\begin{itemize}
			\item For $1\le i,j\le 2$ we have that $Q_{N,ij}^k(r,t)=Q_{N,ji}^k(r,t)$ and $G_{N,ij}^k(r,t)=G_{N,ji}^k(r,t)^t$.
			\item The variables $(Q_{N,nm}^k(r,t))_{n,m=1}^2$, $(V_{N,i1}^k(r,t),V_{N,i2}^k(r,t))$ for $1\le i\le 2$, and $(G_{N,ij}^k(r,t))$ for $1\le i,j\le 2$ are independent.
			\item For each $1\le i,j\le 2$, the entries of $G_{N,ij}^k(r,t)$, except for those required by symmetry, are independent. Moreover they satisfy
			\[\E[[G^k_{N,11}(r,t)]_{ij}^2]=N^{-1}(1+\delta_{ij})\frac{(p-1)}{p},\;\;\E[[G^k_{N,22}(r,t)]_{ij}^2]=N^{-1}(1+\delta_{ij})\frac{(q-1)}{q},\]
			\[\E[G^k_{N,12}(r,t)^2]=N^{-1}.\]
			In particular, $G^k_{N,ii}(r,t)$ is proportional to a $GOE(N_i-2)$ matrix in law.
			\item For $1\le \ell\le 2$, and $1\le i\le N_\ell-2$, the vectors $(V^k_{N,\ell 1}(r,t)_i,V^k_{N,\ell 2}(r,t)_i)$ are i.i.d with covariance given by a fixed matrix only dependant on $(k,p,q,r,t)$.
			\item The random variable $(Q^k_{N,11}(r,t),Q^k_{N,12}(r,t),Q^k_{N,22}(r,t))$ is a Gaussian vector with covariance matrix only dependant $(k,p,q,r,t)$.
		\end{itemize}
	\end{lemma}
	\begin{remark}
		We remark that the correlations between the matrices $M^{1}_{N}(r,t,E_1,E_2)$ and $M^{2}_{N}(r,t,E_1,E_2)$ are not given in the previous lemma, as they will be unimportant in the analysis below.
	\end{remark}
	
	The proof of this result will also be given in Appendix \ref{section:covariance-section}. For $k=1,2$, let $G_N^k:=G_{N}^k(r,t)$ denote the $(N-4)$-by-$(N-4)$ matrix with blocks given by combining the various $G_{N,ij}^k(r,t)$. For $i=1,2$, let us define \[V^i_N=V^i_N(r,t)=\begin{bmatrix}V_{N,11}^i(r,t)&V_{N,12}^i(r,t)\\
		V_{N,21}^i(r,t)&V_{N,22}^i(r,t)\end{bmatrix},\]\[L^i_N=L^i_N(r,t,E_1,E_2)=\begin{bmatrix}M^{i}_N(r,t,E_1,E_2)_{11}&M^{i}_N(r,t,E_1,E_2)_{1N_1}\\M^{i}_N(r,t,E_1,E_2)_{N_11}&M^{i}_N(r,t,E_1,E_2)_{N_1N_1}\end{bmatrix}.\]
	These notations are chosen so that, up to the reordering of the rows and columns of $M^i_N(r,t,E_1,E_2)$ given by $(1,N_1,2,\dots, N-2)$, the matrix $M^i_N(r,t,E_1,E_2)$ is given by
	\[\begin{bmatrix} 
		L^i_N(r,t,E_1,E_2)& V_N^i(r,t)^t\\
		V_N^i(r,t)& G_N^i(r,t)-E_iI\\
	\end{bmatrix}.\label{eqn:ignore-chris}\]
	Here, and in the remainder of this section, we will omit the choice of $(r,t,E_1,E_2)$ from the notation whenever such a choice is clear.
	We will need an upper bound for the absolute determinant of $M_N^i$. To state this, let us denote $v_{r,t}=(\sq{p(1-r^2)},\sq{q(1-t^2)})^t$, and define $W_{N}^i(\epsilon,r,t,E_1,E_2)\ge 0$ by
	\[W^i_N(\epsilon,r,t,E_1,E_2)^2=\frac{4}{\|v_{r,t}\|^2}(\frac{1}{\epsilon^2}\|V^i_N\|_{HS}^4+\|L^i_N\|_{HS}^2)(\frac{1}{\epsilon^2}\|V^i_N v_{r,t}\|^2\|V^i_N\|_{HS}^2+\|L_N^iv_{r,t}\|^2).\]
	Finally, we will use the notational shorthand $G^i_{N,\epsilon}:=G^i_N+i\epsilon  I$.
	
	\begin{lemma}\label{lem:Determinant-Upper-Bound}
		With the notation of Lemma \ref{lem:Hessian-covariance-statement}, for $i=1,2$ and $\epsilon>0$, we have a.s. that \[|\det(M_{N}^i(r,t,E_1,E_2))|\le
		|W_N^i(\epsilon,r,t,E_1,E_2)||\det(G^i_{N,\epsilon}-E_iI)|.\label{determinant-upper bound}\]
	\end{lemma}
	\begin{proof}
		We begin by observing that as for $x\in \R$ we have that $|x|\le |x+i\epsilon|$, we have that
		\[|\det(M_N^i)|\le |\det(M_N^i+i\epsilon I)|.\]
		We note as well that $G_{N,\epsilon}^i$ is invertible, so by the Schur complement formula and (\ref{eqn:ignore-chris}) we have that
		\[|\det(M_N^i+i\epsilon)|=|\det(G^i_{N,\epsilon}-E_i I)||\det(L^i_N-(V^i_N)^t(G_{N,\epsilon}^i)^{-1}V^i_N)|.\]
		If we denote
		\[J^i_{N,\epsilon}:=L_N^i-(V_N^i)^t(G_{N,\epsilon}^i)^{-1}V^i_N,\]
		then it is sufficient to show that
		\[|\det(J^i_{N,\epsilon})|\le |W_N^i(\epsilon)|.\]
		For this, we first observe that for any $2$-by-$2$ matrix, $A$, and any $v,w\in S^1$, with $(v,w)=0$, we have by Hadamard's inequality that
		\[|\det(A)|\le \|Av\|\|Aw\| \le \|A v\|\|A\|_{HS}.\]
		In particular, we have that
		\[|\det(J_{N,\epsilon}^i)|\le \frac{1}{\|v_{r,t}\|}\|J^i_{N,\epsilon}v_{r,t}\|\|J^i_{N,\epsilon}\|_{HS}.\label{eqn:ignore-curt2}\]
		We also observe that for $v,w\in \R^2$, we have that $\|G_{N,\epsilon}^{-1}\|\le \epsilon^{-1}$ so that
		\[|((V^i_N)^tG_{N,\epsilon}^{-1}V^i_N v,w)|\le \epsilon^{-1}\|V^i_Nw\|\|V^i_N v\|.\label{eqn:ignore-curt}\]
		In particular,
		\[\|(V^i_N)^tG_{N,\epsilon}^{-1}V_N^i v\|^2\le \epsilon^{-1}\|V_N^i(V_N^i)^t G_{N,\epsilon}^{-1}V_N^i v\| \|V_N^i v\|\le \]
		\[\epsilon^{-2}\|V_N^i(V_N^i)^t \| \|V_N^i v\|^2\le \epsilon^{-2}\|V_N^i\|_{HS}^2 \|V_N^i v\|^2.\]
		Applying the inequality $(x+y)^2\le 2(x^2+y^2)$, we see that
		\[\|J^i_{N,\epsilon}v_{r,t}\|^2\le 2(\|L^i_{N}v_{r,t}\|^2+\|(V^i_N)^tG_{N,\epsilon}^{-1}V^i_N v_{r,t}\|^2)\le 2(\epsilon^{-2}\|V^i_N v_{r,t}\|^2\|V^i_N\|_{HS}^2+\|L_N^iv_{r,t}\|^2)\]
		Similarly using that $\|A\|_{HS}^2=\sum_{\ell=1}^{2}\|Ae_\ell\|^2$ and (\ref{eqn:ignore-curt}) again
		\[\|J^i_{N,\epsilon}\|_{HS}^2\le 2(\|L^i_N\|^2_{HS}+\|(V^i_N)^tG_{N,\epsilon}^{-1}V^i_N \|^2_{HS})\le 2(\|L^i_N\|^2_{HS}+\epsilon^{-2}\|V^i_N\|^4_{HS}).\]
		Combining these and (\ref{eqn:ignore-curt2}) completes the proof.
	\end{proof}
	
	The law of the matrix $G_{N}^k$ is composed of independent, centered Gaussian random variables, independent of $(r,t,E_1,E_2,k)$. We will need the following result, which allows us to control the right-hand side of (\ref{determinant-upper bound}).
	\begin{lemma}\label{lem:G-det upper bound}
		For $\epsilon>0$ and $\ell\ge 1$
		\[\limsup_{N\to \infty}\left(\sup_{E_1,E_2}\frac{1}{(N-4)}\bigg[\log(\E[\prod_{i=1}^2|\det(G_{N,\epsilon}^i-E_iI)|^\ell])-[\ell \Omega_{\epsilon}(E_1)+\ell \Omega_{\epsilon}(E_2)]\bigg]\right)\le 0,\]
		where here $\Omega_{\epsilon}(E):=\Omega_{p,q,\gamma,\epsilon}(E)=\int \log(|\lam-E+i\epsilon|)\mu_{p,q,\gamma}(d\lam)$.
	\end{lemma}
	\begin{proof}
		Let us denote the law of a single $G^i_N$ as $G_N$, and assume that $N\ge 5$. By H\"{o}lder's inequality, it suffices to show that
		\[\limsup_{N\to \infty}\sup_{E}\frac{1}{(N-4)}[\log(\E[|\det(G_{N,\epsilon}-EI)|^{2\ell}])-2\ell \Omega_{\epsilon}(E)]\le 0.\]
		Let us denote the empirical spectral measure of $G_N$ as $\mu_{G_N}$. We observe that as all entries of $G_N$ are independent Gaussian random variables, and as each entry's variance is bounded by $C/(N-4)$ for some constant fixed constant $C>0$, the matrix $G_N$ considered as a vector satisfies the log-Sobolev inequality with constant $C/(N-4)$, by the Bakry-Emery criterion \cite{Log-Sobolev}. Observing that $\log(|x+i\epsilon|)$ is a Lipschitz function with Lipshitz constant $1/(2\epsilon)$, we see that (see Lemma 2.3.1 of \cite{cupbook}) $(N-4)\int\log(|x-E+i\epsilon|)\mu_{G_N}(dx)$ is a Lipschitz function of the symmetric matrix $G_N$ (endowed with the Hilbert-Schmidt norm) with Lipschitz constant bounded by $\sq{N-4}/\epsilon$. Thus we observe that by Herbst's argument (see section 2.3 of \cite{Log-Sobolev})
		\begin{gather}
			\E[|\det(G_{N,\epsilon}-EI)|^{2\ell}]=\E[\exp(2\ell (N-4)\int\log(|x-E+i\epsilon|)\mu_{G_N}(dx))]\le \\
			\exp\left(2\ell (N-4)\int\log(|x-E+i\epsilon|)\E[\mu_{G_N}](dx)+4C\ell^{2}\epsilon^{-2}\right),
		\end{gather}
		so that
		\[\limsup_{N\to \infty}\sup_{E}\frac{1}{(N-4)}[\log(\E[|\det(G_{N,\epsilon}-EI)|^{2\ell}])-2\ell \Omega_{\epsilon}(E)]\le\]
		\[\limsup_{N\to \infty}\sup_{E}[2\ell \int\log(|x-E+i\epsilon|)\E[\mu_{G_N}](dx)-2\ell \Omega_{\epsilon}(E)]\le \frac{2\ell}{2\epsilon}\W_1(\E[\mu_{G_N}],\mu_{p,q,\gamma}),\]
		where, as above, $\W_1$ denotes the $1$-Wasserstein distance. We now observe that the arguments of Lemma 3.5 in \cite{bipartite}, rescaled to the variances of $G_N$, imply that for some $\kappa>0$, we have that
		\[\W_1(\E[\mu_{G_N}],\mu_{p,q,\gamma})\le N^{-\kappa},\]
		which, when combined with the above observation, completes the proof.
	\end{proof}
	
	\begin{lemma}\label{lem:W_k-norm upper bound}
		Fix the notation of Lemma \ref{lem:Hessian-covariance-statement}. For $r,t\in (-1,1)$, let $(U_1(r,t),U_2(r,t))$ be an (independent) centered Gaussian vector with covariance matrix $\Sigma_U(r,t)$ (defined in (\ref{appendix Sigma_E})). Then for each $\epsilon>0$ and $m\ge 1$ there exists $C>0$ such that for $1\le i\le 2$
		\[\E[W_N^i(\epsilon,r,t,U_1(r,t),U_2(r,t))^{2m}]\le C(2-r^2-t^2)^m.\]
	\end{lemma}
	\begin{proof}
		
		By H\"{o}lder's Inequality, the inequality $(x+y)^l\le 2^l(x^l+y^l)$, and the observation that $\|v_{r,t}\|^2\ge (2-r^2-t^2)$, we see that it suffices to show that
		\[\max(\E[\|L^i_N v_{r,t}\|^{2m}],\E[\| V^i_N v_{r,t}\|^{2m}])\le C(2-r^2-t^2)^{m},\]
		and that
		\[\max(\E[\|L^i_N\|^{2m}_{HS}],\E[\|V^i_N \|^{2m}_{HS}])\le C.\label{eqn:ignore-846}\]
		By employing the inequality $(x+y)^l\le 2^l(x^l+y^l)$ again, we see that to show (\ref{eqn:ignore-846}) it suffices to show that for all fixed choices of $v\in \R^2$, there is $C>0$ such that
		\[\max(\E[\|L^i_N v\|^{2m}],\E[\|V^i_N v\|^{2m}])\le C.\]
		As we will below, we shall denote 
		\[\E_{\D h_N(\n),\D h_N(\n(r,t))}[f]=\E[f|\D h_N(\n)=\D h_N(\n(r,t))=0]\] 
		Let us denote
		\[\Sigma_1(r,t):=\]\[\E_{\D h_N(\n),\D h_N(\n(r,t))}[(\frac{1}{p}E_i E_1h_N(\n),\frac{1}{\sq{pq}}E_i E^1h_N(\n))(\frac{1}{p}E_i E_1h_N(\n),\frac{1}{\sq{pq}}E_i E^1h_N(\n))^t],\]
		\[\Sigma_2(r,t):=\]\[\E_{\D h_N(\n),\D h_N(\n(r,t))}[(\frac{1}{\sq{pq}}E_i E^1h_N(\n),\frac{1}{q}E^i E^1h_N(\n))(\frac{1}{\sq{pq}}E_i E^1h_N(\n),\frac{1}{q}E^i E^1h_N(\n))^t].\]
		Then we see that the components of $V_N^i v$ are all independent centered Gaussian random variables, such that the $k$-th component has variance, \[N^{-1}\sigma_1(v)^2:=N^{-1}v^t\Sigma_1(r,t)v,\]
		for $1\le k\le N_1-1$, and variance
		\[N^{-1}\sigma_2(v)^2:=N^{-1}v^t\Sigma_2(r,t)v,\]
		for $N_1\le k\le N-2$. As conditioning on a random variable can only decrease its variance, we observe that by Lemma \ref{lem:hard-covariance-computation} below, for $i=1,2$, we have that $\sigma_i(v)\le C\|v\|$ for some absolute constant $C$. Moreover by Lemma \ref{lem:Covariance bounds around r,t=1} we have some constant $C$ such that $\sigma_i(v_{r,t})\le C(2-r^2-t^2)$. 
		
		We observe as well that $\|V^i_N v\|^2$ is distributed identically to $\frac{1}{N}(\sigma_1(r)\chi_{N_1-1}^2+\sigma_2(r)\bar{\chi}_{N_2-1}^2)$, where $(\chi_{N_1-1}^2,\bar{\chi}_{N_2-1}^2)$ are independent $\chi^2$-random variables of parameter $N_1-1$ and $N_2-1$, respectively. Recalling that (see \cite{chi-ref})
		\[\E[\chi_{k}^{2m}]=(k-1)(k+1)\dots (k-3+2m),\]
		we see that the required statements for $\E[\|V^1_N v\|^{2m}]$ and $\E [\|V^1_N v_{r,t}\|^{2m}]$ follow from the above bounds on $\sigma_i(v)$. The results for $L^i_N$ follow similarly.
	\end{proof}
	
	With these prerequisites established, we are ready to begin the proof of Proposition \ref{prop:KR-2point}.
	
	\begin{proof}[Proof of Proposition \ref{prop:KR-2point}]
		We note that \[\det(D_N)=N^{(N-2)/4}p^{(N_1-1)/2}q^{(N_2-1)/2}=\]\[\exp(N[\frac{1}{4}\log(N)+\frac{\gamma}{2}\log(p)+\frac{(1-\gamma)}{2}\log(q)]+o(N)).\]
		We also recall that
		\[\vol(S^{N-1})=\frac{2\pi^{N/2}}{\Gamma(N/2)}=\exp(N[-\frac{1}{2}\log(N)+\frac{1}{2}\log(2\pi e)]+o(N)),\]
		so that
		\[C_N:=\exp(N_1[-\log(N_1)+\log(2\pi e)]+N_2[-\log(N_1)+\log(2\pi e)]+o(N))=\]
		\[\exp(N[-\log(N)-\gamma\log(\gamma)-(1-\gamma)\log(1-\gamma)+\log(2\pi e)]+o(N)).\]
		Let us define
		\[g_N(r,t):=g_{N_1,N_2,p,q}(r,t)=\frac{N_1-2}{2}\log(\frac{1-r^2}{1-t^{2q}r^{2p-2}})+\frac{N_2-2}{2}\log(\frac{1-t^2}{1-r^{2p}t^{2q-2}}),\]
		We will further denote the law of \[(D^{-1}_N\D^2h_N(\n)D^{-1}_N,D^{-1}_N\D^2h_N(\n(r,t))D^{-1}_N,h_N(\n),h_N(\n(r,t))),\]
		conditioned on $(\D h_N(\n)=0,\D h_N(\n(r,t))=0)$ by $(H_{N,1}^D(r,t),H_{N,2}^D(r,t),U_1(r,t),U_2(r,t))$.
		
		Combining the above computations with Lemma \ref{lem:Gradient-covariance-statement} and Lemma \ref{lem:Hessian-covariance-statement}, we obtain that
		\[\E[\Crit_{N,2}(B,I_O)]=\bar{C}_N\int_{I_O}\exp(g_N(r,t))\times \]\[f_L(r,t)((1-r^2)(1-t^2))^{-1/2} \E[\prod_{i=1}^2|\det(H_{N,i}^D(r,t))|I(U_i(r,t)\in \sq{N}B)]drdt \label{proof-Crit-formula},\]
		where here $\bar{C}_N=\exp(2NC_{p,q,\gamma}+o(N))$.
		
		Choose $\epsilon>0$ and $m\ge 4$ and set $\ell=\ell(m):=m/(m-1)$. In what follows, we will work with constants $C>0$, which we will allow to depend on the choice of $\epsilon$ and $m$, but not $N$. In addition, the choice of these constants will be allowed to change line by line. Now we have by Lemma \ref{lem:Determinant-Upper-Bound} and H\"{o}lder's inequality that
		\[\E[\prod_{i=1}^2|\det(H_{N,i}^D(r,t))|I(U_i(r,t)\in \sq{N}B)]\le (\cal{E}^{(1)}_{\epsilon}(r,t))^{1/\ell}(\cal{E}_{\epsilon}^{(2)}(r,t))^{1/2m},\label{eqn:hold-proof}\]
		where here
		\[\cal{E}_{\epsilon}^{(1)}(r,t)=\E[\prod_{i=1}^2|\det(G^i_{N,\epsilon}-U_i(r,t)/\sq{N})|^\ell I(U_i(r,t)/\sq{N}\in B)],\]
		\[\cal{E}^{(2)}_{\epsilon}(r,t)=\prod_{i=1}^2\E[W_N^i(\epsilon,r,t,U_1(r,t),U_2(r,t))^{2m}],\]
		and where $G^i_{N,\epsilon}$ is independent of $U_i(r,t)$. Applying Lemma \ref{lem:W_k-norm upper bound} we see that
		\[(\cal{E}^{(2)}_{\epsilon}(r,t))^{1/2m}\le C(2-r^2-t^2).\]
		We recall from Lemma \ref{lem:upper bound on H_1,H^1 covariance} that
		\[f_{L}(r,t)\le C(2-r^2-t^2)^{-1}.\]
		Combining this with (\ref{eqn:hold-proof}), we see that
		\[\E[\Crit_{N,2}(B,I_O)]\le C \bar{C}_N\int_{I_O} \frac{\exp(g_N(r,t))\cal{E}_{\epsilon}^{(1)}(r,t)^{1/\ell}}{\sq{(1-r^2)(1-t^2)}}drdt.\]
		An additional application of H\"{o}lder's inequality shows that
		\[\int_{I_O} \frac{\exp(g_N(r,t))\cal{E}_{\epsilon}^{(1)}(r,t)^{1/\ell}}{\sq{(1-r^2)(1-t^2)}}drdt\le C(\int_{I_O} \frac{\exp(\ell g_N(r,t))\cal{E}_{\epsilon}^{(1)}(r,t)}{\sq{(1-r^2)(1-t^2)}}drdt)^{1/\ell}.\]
		In view of the uniform bound of Lemma \ref{lem:G-det upper bound} we see that for $N$ sufficiently large (independent of $(r,t)$) we have that
		\[\cal{E}_{\epsilon}^{(1)}(r,t)\le e^{N\epsilon}\E[\prod_{i=1}^2\exp(\ell(N-4)\Omega_{\epsilon}(U_i(r,t)/\sq{N}))I(U_i(r,t)/\sq{N}\in B)].\]
		We now note that
		\[|\Omega_{\epsilon}(E/\sq{N})-\Omega_{\epsilon}(E/\sq{N-4})|\le \int^{E/\sq{N-4}}_{E/\sq{N}}|\Omega'_{\epsilon}(y)|dy.\]
		We also note as well that
		\[|\Omega'_{\epsilon}(y)|\le \int \frac{|y+\lam|}{|y+\lam+i\epsilon|^2}\mu_{p,q,\gamma}(d\lam).\]
		Noting that $\mu_{p,q,\gamma}$ has compact support, we easily see that for fixed $\epsilon>0$ we may choose $C$ so that for $y\in \R$
		\[|\Omega'_{\epsilon}(y)|\le C\frac{1}{|y|+1},\]
		and in particular, we see that
		\[\sup_{E\in \R}|\Omega_{\epsilon}(E/\sq{N})-\Omega_{\epsilon}(E/\sq{N-4})|\le C/N.\]
		In addition, we see that for sufficiently large $N$ we have that $ (N/(N-4))^{1/2}B\subseteq B_{\epsilon}(B)$, where here $B_{\epsilon}(B)=\{x\in \R:d(x,B)<\epsilon\}$. Thus for large enough $N$ (independent of $(r,t)$) we have that
		\[\cal{E}_{\epsilon}^{(1)}(r,t)\le Ce^{N\epsilon}\E[\prod_{i=1}^2\exp(\ell (N-4)\Omega_{\epsilon}(U_i(r,t)/\sq{N-4}))I(U_i(r,t)/\sq{N-4}\in B_{\epsilon}(B))].\]
		Recalling the definition of $g_{p,q,\gamma}$ given in (\ref{appendix Sigma_E}), and denoting $g=g_{p,q,\gamma}$ we see that
		\[g_N(r,t)-(N-4)(1-\epsilon)g(r,t)=C_{N,1,\epsilon}\log(\frac{1-r^2}{1-t^{2q}r^{2p-2}})+C_{N,2,\epsilon}\log(\frac{1-t^2}{1-r^{2p}t^{2q-2}}),\]
		\[C_{N,1,\epsilon}:=\frac{N_1-2-(N-4)(1-\epsilon)\gamma}{2},\;\;\; C_{N,2,\epsilon}:=\frac{N_2-2-(N-4)(1-\gamma)(1-\epsilon)}{2}.\]
		Both $\log(\frac{1-r^2}{1-t^{2q}r^{2p-2}})$ and $\log(\frac{1-t^2}{1-r^{2p}t^{2q-2}})$ are bounded above. In addition, for each $k=1,2$, we have that $ N \epsilon\ge C_{N,k,\epsilon}\ge 0$ for large enough $N$. From this, we see that for sufficiently large $N$ we have that 
		\[g_N(r,t)\le (N-4)(1-\epsilon)g(r,t)+2C_0N\epsilon,\]
		where here
		\[C_0=\sup_{r,t\in (-1,1)}|\log(\frac{1-r^2}{1-t^{2q}r^{2p-2}})|+\sup_{r,t\in (-1,1)}|\log(\frac{1-t^2}{1-r^{2p}t^{2q-2}})|<\infty.\]
		Combining these observations, we see that there is $C$ such that for large enough $N$
		\[\int_{I_O} \frac{\exp(\ell g_N(r,t))\cal{E}_{\epsilon}^{(1)}(r,t)}{\sq{(1-r^2)(1-t^2)}}drdt\le e^{2N(C_0+1)\epsilon}\int_{I_O}\frac{\exp(\ell (N-4)(1-\epsilon)g(r,t))}{\sq{(1-r^2)(1-t^2)}}\times \]\[\E[\prod_{i=1}^2\exp(\ell(N-4)\Omega_{\epsilon}(U_i(r,t)/\sq{N-4}))I(U_i(r,t)/\sq{N-4}\in B)]drdt.\label{equation:in proof of Varadhan}\]
		We will analyze the latter integral through an application of Varadhan's Lemma (see Theorem 4.3.1 and Lemma 4.3.6 of \cite{Varadhan}). Let us fix a pair of independent standard Gaussian random variables $(X_1,X_2)$, and observe that
		\[(X_1,X_2)\Sigma_U(r,t)^{1/2}\stackrel{d}{=}(U_1(r,t),U_2(r,t)).\]
		Furthermore, let us define a random vector $(R,T)$ on $I_O$, independent of $(X_1,X_2)$, with pdf proportional to $(1-r^2)^{-1/2}(1-t^2)^{-1/2}$. Let us additionally define
		\[W=\{(r,t,x_1,x_2)\in I_O\times \R^2:(x_1,x_2)\Sigma_U(r,t)^{1/2}\in B_{\epsilon}(B)^2\}.\]
		Let us also define \[\phi_{\epsilon}(r,t,E_1,E_2)=(1-\epsilon)g(r,t)+\sum_{i=1}^2\Omega_{\epsilon}([(E_1,E_2)\Sigma_U(r,t)^{1/2}]_i).\]
		We observe that up to a fixed constant, we may rewrite the right-hand side of (\ref{equation:in proof of Varadhan}) as
		\[e^{2N(C_0+1)\epsilon}\E[\exp(\ell(N-4)\phi_{\epsilon}(R,T,\frac{X_1}{\sq{N-4}},\frac{X_2}{\sq{N-4}}))I((R,T,\frac{X_1}{\sq{N-4}},\frac{X_2}{\sq{N-4}})\in W)].\]
		We now observe that the random vector $(R,T,X_1/\sq{N},X_2/\sq{N})$ satisfies a LDP on $\bar{I}_O\times \R^2$ with good rate function \[J(r,t,u_1,u_2)=\frac{u_1^2}{2}+\frac{u_2^2}{2}.\]
		To apply Varadhan's Lemma, we will need to first verify a moment condition. Let $\ell'>0$ be arbitrary. To begin, we again apply H\"{o}lder's inequality to see that
		\[\E[\prod_{i=1}^2\exp(\ell'N\Omega_{\epsilon}(U_i(r,t)/\sq{N}))]\le \prod_{i=1}^2\E[\exp(2\ell'N\Omega_{\epsilon}(U_i(r,t)/\sq{N}))]^{1/2}.\label{equation:proof partial bound moment}\]
		We note that as the support of $\mu_{p,q,\gamma}$ is compact, there $C>0$ such that for $x\in \R$
		\[\Omega_{\epsilon}(x)\le C\log(|x+i|)+C.\]
		As the variance of a conditioned Gaussian random variable is bounded by the variance of the original, we have that $\E[U_i(r,t)^2]\le 1$ for $i=1,2$. Noting as well that for $a\ge 1$, we have that $\log(|x+i|)\le \log(|ax+i|)$, we see that
		\[\E[\exp(2\ell'N\Omega_{\epsilon}(U_i(r,t)/\sq{N}))]\le e^{2\ell' NC}\E[\exp(C2\ell'N\log(|X_i/\sq{N}+i|))]=\]
		\[\sq{\frac{N}{2\pi}}\int  \exp(N[2\ell'C+ 2\ell'C\log(|x+i|)-x^2/2])dx.\label{equation:proof moment bound}\]
		For large enough $C$, we see that $2\ell'C\log(|x+i|)\le C^2+x^{2}/4$, and plugging this into the right-hand side of (\ref{equation:proof moment bound}) we see that
		\[\limsup_{N\to\infty}\frac{1}{N}\log(\E[\exp(2\ell'N\Omega_{\epsilon}(U_i(r,t)/\sq{N}))])\le 2\ell' C+C^2 <\infty.\]
		Combining this with (\ref{equation:proof partial bound moment}) and the above observation that $g(r,t)$ is bounded above, we may easily obtain that
		\[\limsup_{N\to \infty}\frac{1}{N}\log(\E[\exp(\ell'N\phi_{\epsilon}(R,T,\frac{X_1}{\sq{N}},\frac{X_2}{\sq{N}}))])<\infty.\]
		We now note that $\phi_{\epsilon}$ is continuous on $(-1,1)^2\times \R^2$. The standard form of Varadhan's Lemma we wish to apply (see Theorem 4.3.8 of \cite{Varadhan}) would require that $\phi_{\epsilon}$ admits a continuous extension to $[-1,1]^2 \times \R^2$ (or at least $\bar{I}_O\times \R^2$).
		
		To get around this we note that with the above bounds on $\E[U_i(r,t)^2]$ and $g$, we have that \[\sup_{r,t\in (-1,1)}\phi_{\epsilon}(r,t,E_1,E_2)<\infty.\]
		Thus $\phi_{\epsilon}$ may be extended to an upper semi-continuous function on $\bar{\phi}_{\epsilon}:[-1,1]^2\times \R^2\to \R$, such that $\sup_{r,t\in (-1,1)}\phi_{\epsilon}(r,t,E_1,E_2)=\sup_{r,t\in (-1,1)}\bar{\phi}_{\epsilon}(r,t,E_1,E_2)$ for each $E_1,E_2\in\R$. Then applying Lemmas 4.3.6 and 4.3.8 of \cite{Varadhan} we see that
		\[\limsup_{N\to \infty}\frac{1}{N}\log(\E[\exp(\ell N\phi_{\epsilon}(R,T,\frac{X_1}{\sq{N}},\frac{X_2}{\sq{N}}))I((R,T,\frac{X_1}{\sq{N}},\frac{X_2}{\sq{N}})\in W)])\le \]
		\[\sup_{(r,t,x_1,x_2)\in W}[2(C_0+1)\epsilon+\ell \phi_{\epsilon}(r,t,x_1,x_2)-\frac{x_1^2}{2}-\frac{x_2^2}{2}]=2(C_0+1)\epsilon+\]
		\[\sup_{E_1,E_2\in B_{\epsilon}(B),r,t\in I_O,}[\ell(1-\epsilon) g(r,t)+\ell\Omega_{\epsilon}(E_1)+\ell\Omega_{\epsilon}(E_2)-\frac{1}{2}(E_1,E_2)\Sigma_U(r,t)^{-1}(E_1,E_2)^t].\]
		Combining this with the above bounds, we see that
		\begin{gather}
			\limsup_{N\to \infty}\frac{1}{N}\log(\E[ \Crit_{N,2}(B,I_O)])\le\\
			C_{p,q,\gamma}+2(C_0+1)\epsilon+\sup_{E_1,E_2\in B_{\epsilon}(B),r,t\in I_O,}[C\epsilon+\ell (1-\epsilon)g(r,t)+\\
			\ell\Omega_{\epsilon}(E_1)+\ell\Omega_{\epsilon}(E_2)-\frac{1}{2}(E_1,E_2)\Sigma_U(r,t)^{-1}(E_1,E_2)^t].\label{equation:proof epsilon-l versions of Psi}
		\end{gather}
		Observing that $\Sigma_U(r,t)\le I$ and that $\Omega_{\epsilon}(x)=O(\log(|x+i\epsilon|))$, one easily sees that taking $m\to \infty$ (so that $\ell\to 1$) and then $\epsilon\to 0$, the right-hand side of (\ref{equation:proof epsilon-l versions of Psi}) becomes 
		\[\sup_{r,t\in I_O}\sup_{E_1,E_2\in B}\Sigma_{p,q,\gamma,2}(r,t,E_1,E_2),\]
		which is the desired upper bound.
	\end{proof}
	
	\appendix
	
	\section{Covariance Computations and Conditional Densities \label{section:covariance-section}}
	
	In this appendix we will study the covariance structure of the Gaussian vector \[(h_N(\n),\D h_N(\n), \D^2 h_N(\n),h_N(\n(r,t)),\D h_N(\n(r,t))),\]
	after which we will investigate some of its conditional density functions.
	
	We will be interested specifically in studying the case where the orthonormal frame field $(E_i)_{i=1}^{N-2}$ is induced by a pullback of a choice of an orthonormal frame field on each component of the product. That is, let us write $E^i=E_{i+N_1-1}$, and assume that $(E_i)_{i=1}^{N_1-1}$ is obtained from an orthonormal frame field on the factor $S^{N_1-1}$, and similarly that $(E^i)_{i=1}^{N_2-1}$ is obtained from an orthogonal frame field on the factor $S^{N_2-1}$. We observe that in this case, the action of $E^i$ commutes with $E_j$ for each $i,j$. 
	
	With $\delta_{ij}$ denoting the standard Kronecker $\delta$, we denote $\delta_{ijk}=\delta_{ij}\delta_{jk}$, $\delta_{ij\neq 1}=\delta_{ij}(1-\delta_{i1})$, etc. In addition, we employ the shorthand $x_*=\sq{1-x^2}$.
	\begin{lemma}\label{lem:hard-covariance-computation}
		For $p,q\ge 2$,  $0<\gamma<1$ and $r,t\in [-1,1]$ there exists an orthogonal frame field $E$, such that:
		\[\E[h_N(\n)h_N(\n(r,t))]=r^p t^q\]
		\[\E[E_i h_N(\n)h_N(\n(r,t))]=-\E[h_N(\n)E_i h_N(\n(r,t))]=pr^{p-1}r_*t^{q}\delta_{1i}\]
		\[\E[E^i h_N(\n)h_N(\n(r,t))]=-\E[h_N(\n)E^i h_N(\n(r,t))]=qt^{q-1}t_*r^p\delta_{1i}\]
		\[\E[E_i h_N(\n)E_j h_N(\n(r,t))]=t^q(pr^{p-1}\delta_{ij\neq 1}+\delta_{ij1}pr^{p-2}(pr^2-(p-1)))\]
		\[\E[E^i h_N(\n)E^j h_N(\n(r,t))]=r^p(qt^{q-1}\delta_{ij\neq 1}+\delta_{ij1}qt^{q-2}(qt^2-(q-1)))\]
		\[\E[E^i h_N(\n)E_j h_N(\n(r,t))]=\E[E_i h_N(\n)E^j h_N(\n(r,t))]=-pqt^{q-1}r^{p-1}r_*t_*\delta_{ij1}\]
		\[\E[E_i E_j h_N(\n)h_N(\n(r,t))]=\E[h_N(\n)E_i E_j h_N(\n(r,t))]=t^q(-pr^p\delta_{ij}+p(p-1)r^{p-2}r_*^2\delta_{ij1})\]
		\[\E[E^i E^j h_N(\n)h_N(\n(r,t))]=\E[h_N(\n)E^i E^j h_N(\n(r,t))]=r^p(-qt^q\delta_{ij}+q(q-1)t^{q-2}t_*^2\delta_{ij1})\]
		\[\E[E^i E_j h_N(\n)h_N(\n(r,t))]=\E[h_N(\n)E^i E_j(\n(r,t))]=pqr^{p-1}t^{q-1}r_*t_*\delta_{ij1}\]
		\[\E[E^i E^j h_N(\n)E_k h_N(\n(r,t))]=-\E[E_k h_N(\n)E^i E^j h_N(\n(r,t))]=\]\[-pr_*r^{p-1}\delta_{k1}(-\delta_{ij}t^q q+\delta_{ij1}q(q-1)t^{q-2}t_*^2)\]
		\[\E[E_i E_j h_N(\n)E^k h_N(\n(r,t))]=-\E[E^k h_N(\n)E_i E_j h_N(\n(r,t))]=\]\[-qt_*t^{q-1}\delta_{k1}(-\delta_{ij}r^p p+\delta_{ij1}p(p-1)r^{p-2}r_*^2)\]
		\[\E[E_i E^j h_N(\n)E_k h_N(\n(r,t))]=-\E[E_k h_N(\n)E_i E^j h_N(\n(r,t))]=\]\[q\delta_{j1}t^{q-1}t_*(pr^{p-1}\delta_{ik\neq 1}+\delta_{ik1}pr^{p-2}(pr^2-(p-1)))\]
		\[\E[E^i E_j h_N(\n)E^k h_N(\n(r,t))]=-\E[E^k h_N(\n)E^i E_j h_N(\n(r,t))]=\]\[p\delta_{j1}r^{p-1}r_*(qt^{q-1}\delta_{ik\neq 1}+\delta_{ik1}qt^{q-2}(qt^2-(q-1)))\]
		\[\E[E_i E_j h_N(\n)E_k h_N(\n(r,t))]=-\E[E_k h_N(\n)E_i E_j h_N(\n(r,t))]=\]\[t^q[\delta_{k\neq 1}(\delta_{i1}\delta_{jk}+\delta_{j1}\delta_{ik})p(p-1)r_*r^{p-2}+\delta_{k1}(\delta_{ij}p^2r^{p-1}r_*+\delta_{ij1}p(p-1)r^{p-3}r_*(2r^2-(p-2)r^2_*))]\]
		\[\E[E^i E^j h_N(\n)E^k h_N(\n(r,t))]=-\E[E^k h_N(\n)E^i E^j h_N(\n(r,t))]=\]\[r^p[\delta_{k\neq 1}(\delta_{i1}\delta_{jk}+\delta_{j1}\delta_{ik})q(q-1)t_*t^{q-2}+\delta_{k1}(\delta_{ij}q^2t^{q-1}t_*+\delta_{ij1}q(q-1)t^{q-3}t_*(2t^2-(q-2)t^2_*))]\]
	\end{lemma}
	
	\begin{proof}
		The proof will proceed similarly to Lemma 30 of \cite{pspin-second}. Let us denote $r_1=r$ and $r_2=t$. Let $P_{\n,N}:S^{N-1}\to \R^{N-1}$ denote projection away from the first coordinate, let $\theta_i\in [-\pi/2,\pi/2]$ denote the angles such that $\sin(\theta_i)=r_i$, and let $R_{\varphi,N}:S^{N-1}\to S^{N-1}$ denote the rotation mapping
		\[R_{\varphi,N}(x_1,\dots, x_N)=(\cos(\varphi)x_1+\sin(\varphi)x_2,-\sin(\varphi)x_1+\cos(\varphi)x_2,x_3,\dots, x_N).\]
		
		For $i=1,2$, take a neighborhood $U_i$ of $\n_{N_i}$ in $S^{N_i-1}$ such that restriction of $P_{\n,N_i}$ is a diffeomorphism onto its image. Similarly take a neighborhood $V_i$ of $\n_{N_i}(r_i)$ in $S^{N_i-1}$ such that $P_{\n,N_i}\circ R_{-\theta_i,N_i}$ is a diffeomorphisms onto its image. We denote these images as $\bar{U}_i$ and $\bar{V}_i$, respectively. We define functions $\bar{h}_1:\bar{U}_1\times \bar{U}_2\to \R$ and $\bar{h}_2:\bar{V}_1\times \bar{V}_2\to \R$ by
		\[\bar{h}_{1}=h_N\circ (P_{\n,N_1}\times P_{\n,N_2})^{-1};\;\;\;\bar{h}_{2}=h_N\circ (P_{\n,N_1}\circ R_{-\theta_1,N_1}\times P_{\n,N_2}\circ R_{-\theta_2,N_2})^{-1}.\]
		
		In the proof of Lemma 30 of \cite{pspin-second} (and more specifically in their ft. 5) they show that for each $N\ge 2$ and $\theta\in [-\pi/2, \pi/2]$, there is an orthonormal frame field, $E^{N,\theta}=(E_i^{N,\theta})_{i=1}^{N-1}$, on $S^{N-1}$ with the following property: for any smooth function $f:S^{N-1}\to \R$ if we denote
		\[\sin(\theta)=r,\;\;\;\bar{f}_1=f\circ P_{\n,N}^{-1},\;\;\; \bar{f}_2=f\circ (P_{\n,N}\circ R_{-\theta,N})^{-1},\]
		then for $1\le i,j\le N-1$ we have that
		\[E^{N,\theta}_i f(\n)=\partial_i\bar{f}_1(0),\;\;E^{N,\theta}_i E^{N,\theta}_j f(\n)=\partial_i \partial_j \bar{f}_1(0),\label{eqn:frame-field-1}\]
		\[E^{N,\theta}_i f(\n(r))=\partial_i \bar{f}_2(0),\;\;E^{N,\theta}_i E^{N,\theta}_j f(\n(r))=\partial_i \partial_j \bar{f}_2(0),\label{eqn:frame-field-2}\]
		where $\partial_i$ denotes the $i$-th standard Euclidean partial derivative. We define our orthogonal frame field $E$ on $S^{N_1-1}\times S^{N_2-1}$ by letting $(E_i)_{i=1}^{N_1-1}=(E_i^{N_1,\theta_1})_{i=1}^{N_1-1}$, where $E_i^{N_1,\theta_1}$ acts on the coordinates of the first sphere in the product, and $(E^i)_{i=1}^{N_2-1}=(E_i^{N_2,\theta_2})_{i=1}^{N_2-1}$, where $E_i^{N_2,\theta_2}$ acts on the coordinates of the second sphere. Let us write the $i$-th Euclidean coordinate for $\bar{U}_1\subseteq \R^{N_1-1}$ as $x_i$ and the $i$-th Euclidean coordinate for $\bar{U}_2\subseteq \R^{N_2-1}$ as $y_i$. Then employing (\ref{eqn:frame-field-1}) and (\ref{eqn:frame-field-2}) we see that for $1\le i,j\le N_1-1$ and $1\le k,l\le N_2-1$
		\[
		(h_N(\n),E_i h_N(\n), E^k h_N(\n), E_i E_j h_N(\n), E_i E^k h_N(\n), E^k E^l h_N(\n))=
		\]
		\[
		(\bar{h}_{1}(0), \frac{d}{d x_i}\bar{h}_{1}(0), \frac{d}{d y_k}\bar{h}_1(0), \frac{d^2}{d x_i d x_j}\bar{h}_{1}(0), \frac{d^2}{d y_k d x_i}\bar{h}_1(0),  \frac{d^2}{d y_k d y_l} \bar{h}_1(0)).
		\]
		A similar relationship holds for between the derivatives of $h_N$ at $\n(r,t)$ and $\bar{h}_2$ at $0$.
		
		If we now define
		\[\rho_{r,N}(x,y)=\sum_{i=2}^{N-1}x_iy_i+rx_{1}y_{1}+r_*x_{1}\|y\|_*-r_*y_{1}\|x\|_*+r\|x\|_*\|y\|_*,\]
		and let $C(x,y,z,w):=\E[\bar{h}_1(x,y)\bar{h}_2(z,w)]$, then for $x,z\in\R^{N_1-1}$ and $y,w\in \R^{N_2-1}$, we have
		\[C(x,y,z,w)=(\rho_{r,N_1}(x,z))^p(\rho_{t,N_2}(y,w))^q.\]
		
		We recall (5.5.4) of \cite{AT}, which states that for an arbitrary centered Gaussian field on $\R^n$,  $h$, with smooth covariance function, $C$, we have that
		\[\E[\frac{d^k}{dx_{i_1}\dots dx_{i_k}}h(x)\frac{d^l}{dy_{j_1}\dots dy_{j_l}}h(y)]=\frac{d^k}{dx_{i_1}\dots dx_{i_k}}\frac{d^l}{dy_{j_1}\dots dy_{j_l}}C(x,y).\label{covariance of derivative formula}\]
		Routine calculation using these formulas now yields the desired results.
	\end{proof}
	
	From this result we observe the following corollary.
	
	\begin{corollary}
		\label{corr:ind cor}
		Let us take $r,t\in [-1,1]$ and the choice of frame field in Lemma \ref{lem:hard-covariance-computation}. Then all the following random variables are independent.
		\begin{itemize}
			\item $E_i E_j h_N(\n)+p\delta_{ij}h_N(\n)$ for $1<i\le j\le N_1-1$. 
			\item $E^i E^j h_N(\n)+q\delta_{ij}h_N(\n)$ for $1<i\le j\le N_2-1$.
			\item $E_i E^j h_N(\n)$ for $1<i\le N_1-1,1<j\le N_2-1$.
			\item $(E_{1}E_i h_N(\n),E^1 E_i h_N(\n),E_i h_N(\n),E_i h_N(\n(r,t)))$ for any $1<i\le N_1-1$.
			\item $(E_{1}E^i h_N(\n),E^1 E^i h_N(\n),E^i h_N(\n),E^i h_N(\n(r,t)))$ for any $1<i\le N_2-1$.
			\item $(E_{1}E_1h_N(\n),E_{1}E^1h_N(\n),E^{1}E^1h_N(\n),h_N(\n),h_N(\n(r,t))$,\\
			$E_1h_N(\n),E^1h_N(\n),E_1h_N(\n(r,t)),E^1h_N(\n(r,t)))$.
		\end{itemize}
	\end{corollary}
	
	We will now turn our attention to the computation of the relevant conditional densities, beginning with that of $(\D h_N(\n),\D h_N(\n(r,t)))$. In what follows we will repeatedly use the formula for the conditional law of Gaussian distributions (see for example, (1.2.7) and (1.2.8) of \cite{AT}). We observe that for $1<i<N_1-1$ and $1<j<N_2-1$
	\[\E[(E_i h_N(\n), E_i h_N(\n(r,t)))(E_i h_N(\n), E_i h_N(\n(r,t)))^t]=\begin{bmatrix}
		p & p r^{p-1} t^q\\
		p r^{p-1} t^q & p
	\end{bmatrix},\label{eqn:off-1}\]
	\[\E[(E^j h_N(\n), E^j h_N(\n(r,t)))(E^j h_N(\n), E^j h_N(\n(r,t)))^t]=\begin{bmatrix}
		q & q r^{p} t^{q-1}\\
		q r^{p} t^{q-1} & q
	\end{bmatrix}.\label{eqn:off-2}\]
	From these we see that the vectors $(E_i h_N(\n), E_i h_N(\n(r,t)))$ and $(E^j h_N(\n), E^j h_N(\n(r,t)))$ are non-degenerate unless $|r|=|t|=1$. The remaining entries are more complicated. We may compute that \[\E[(E_1h_N(\n),E^1h_N(\n))(E_1h_N(\n(r,t)),E^1h_N(\n(r,t)))^t]=\]
	\[\begin{bmatrix}
		t^q p  r^{p-2}(p r^2-(p-1))& -pqr^{p-1}t^{q-1}r_* t_*\\
		-pqr^{p-1}t^{q-1}r_* t_* & r^p q t^{q-2}(q t^2-(q-1))
	\end{bmatrix}=:\bar{\Sigma}_{L}(r,t),
	\]
	so that if we denote the covariance matrix of the vector
	\[(E_1h_N(\n),E^1h_N(\n),E_1h_N(\n(r,t)),E^1h_N(\n(r,t)))\label{eqn:ignore-7}\]
	as $\Sigma_L(r,t)$ then we have that
	\[
	\Sigma_L(r,t):=
	\begin{bmatrix}
		\bar{\Sigma}_{L}(1,1)& \bar{\Sigma}_{L}(r,t)\\
		\bar{\Sigma}_{L}(r,t)& \bar{\Sigma}_{L}(1,1)
	\end{bmatrix}.\]
	We note for convenience that
	\[\bar{\Sigma}_{L}(1,1)=\begin{bmatrix}p & 0 \\ 0 & q \end{bmatrix}.\]
	We observe that when $\Sigma_L(r,t)$ is invertible, the density of the vector (\ref{eqn:ignore-7}) at zero is
	\[\varphi_{E_1h_N(\n),E^1h_N(\n),E_1h_N(\n(r,t)),E^1h_N(\n(r,t))}(0,0,0,0)=\frac{1}{(2\pi)^2\sq{\det(\Sigma_{L}(r,t))}}=:f_L(r,t).\label{appendix f_L}\]
	\begin{proof}[Proof of Lemma \ref{lem:Gradient-covariance-statement} and Lemma \ref{lem:Hessian-covariance-statement}]
		By Lemma \ref{lem:Appendix-non-degeneracy} below, for any $r,t\in (-1,1)$ we have that the vector 
		\[(\D h_N(\n),\D h_N(\n(r,t)))\]
		is non-degenerate. Equation (\ref{eqn:graident-to-find}) now follows from (\ref{eqn:off-1}) and (\ref{eqn:off-2})  and the formula for the conditional density of a Gaussian random variable. We will use the notation 
		\[\E_{\D h_N(\n),\D h_N(\n(r,t))}[f]=\E[f|\D h_N(\n)=\D h_N(\n(r,t))=0].\]
		We define
		\[\Sigma_U(r,t):=\E_{\D h_N(\n),\D h_N(\n(r,t))}[(h_N(\n),h_N(\n(r,t)))(h_N(\n),h_N(\n(r,t)))^t].\label{appendix Sigma_E}\]
		Then Lemma \ref{lem:Appendix-non-degeneracy} guarantees that $\Sigma_U(r,t)$ is strictly positive-definite, and the remaining claims follow from Lemma \ref{lem:Appendix-non-degeneracy} and the formulas for the conditional law of Gaussian distributions. Together these results complete the proof of Lemma \ref{lem:Gradient-covariance-statement}.
		
		The proof of Lemma \ref{lem:Hessian-covariance-statement} follows from Lemma \ref{lem:hard-covariance-computation} and Corollary \ref{corr:ind cor} and the same conditional formulas.
	\end{proof}
	The next result shows that the only points where $\det(\Sigma_L(r,t))=0$ have $|r|=|t|=1$.
	
	\begin{lemma}\label{lem:upper bound on H_1,H^1 covariance}
		For $p,q\ge 5$, there is a constant $C>0$, such that for $r,t\in [-1,1]$, we have that
		\[f_{L}(r,t)\le C(2-r^2-t^2)^{-1}.\]
	\end{lemma}
	\begin{proof}
		Appealing to Lemma \ref{lem:Appendix-non-degeneracy} as before, for $r,t\in (-1,1)$ we have that the vector \[(\D h_N(\n),\D h_N(\n(r,t)))\] is non-degenerate, so that in particular $\det(\Sigma_L(r,t))>0$. From this and continuity, we also see that $\det(\Sigma_L(r,t))\ge 0$, so that, in particular, it suffices to show that there is $c>0$ such that
		\[\det(\Sigma_L(r,t))\ge c(2-r^2-t^2)^2.\]
		We observe that
		\[\bar{\Sigma}_{L}(r,\pm 1)=(\pm 1)^{q}\begin{bmatrix}
			pr^{p-2}(pr^2-(p-1))& 0\\
			0 & r^{p}
		\end{bmatrix}.\]
		We note that for $|r|<1$, we have by the AM-GM inequality that
		\[\frac{1}{(p-1)}\frac{1-r^{2p-4}}{1-r^2}=\frac{1}{(p-1)}\sum_{i=0}^{p-2}r^{2i}> |r|^{p-2},\]
		so in particular, $(p-1)r^{p-2}(1-r^2)<1$ and so $1+r^{p-2}(pr^2-(p-1))\ge 1-(p-1)r^{p-2}(1-r^2)>0$. As $|r|<1$, we more easily see that $1-r^{p-2}(pr^2-(p-1))\ge 1-r^{p-2}>0$ and that $|r|^p<1$. In particular, we see that if $|r|<1$ then \[\det(\Sigma_L(r,\pm 1))=(1-pr^{p-2}(pr^2-(p-1)))(1+pr^{p-2}(pr^2-(p-1)))(1-r^{2p})>0\]
		Similarly we see that if $|t|<1$ then $\det(\Sigma_L(\pm 1,t)>0$. Thus what remains is to study the order of vanishing of $\det(\Sigma_L(r,t))$ around the points $|r|=|t|=1$. We observe that $\det(\Sigma_L(r,t))$ is even in both $r$ and $t$, so with the prior positivity results, it is sufficient to show that
		\[\liminf_{(r,t)\to 1}\frac{\det(\Sigma_L(r,t))}{(2-r^2-t^2)^2}>0\]
		where the limit here (and below) is taken with $r,t\in (-1,1)$. Let us denote $\delta=(2-r^2-t^2)$. Then for $r,t\in (-1,1)$ sufficiently close to $(1,1)$, we note that
		\[\bar{\Sigma}_L(r,t)=O(\delta^2)+\begin{bmatrix}
			p& 0\\
			0& q\\
		\end{bmatrix}+\]\[\begin{bmatrix}
			pq(t-1)+p(p^2-(p-1)(p-2))(r-1)& -pq r_*t_*\\
			-pq r_*t_*& pq(r-1)+q(q^2-(q-1)(q-2))(t-1)
		\end{bmatrix},\]
		so that we have that
		\[\det(\Sigma_L(1,1))=0,\D \det(\Sigma_L(1,1))=0.\]
		Finally a long but direct computation yields that
		\[\D^2 \det(\Sigma_L(1,1))=p^2q^2\begin{bmatrix}
			4p(6p-4)& 12(p-1)(q-1)-4\\
			12(p-1)(q-1)-4& 4q(6q-4)\end{bmatrix}.\]
		Let us denote the matrix on the right-hand side as $H_L(p,q)$. We see that by Taylor's theorem
		\[\liminf_{(r,t)\to 1}\frac{\det(\Sigma_L(r,t))}{(2-r^2-t^2)^2}=\liminf_{(r,t)\to 1}\frac{(1-r,1-t)H_L(p,q)(1-r,1-t)^t}{2(2-r^2-t^2)^2}.\]
		On the other hand for $p,q\ge 2$ all entries of $H_L(p,q)$ are positive, so for small enough $c>0$ and $r,t\in (-1,1)$
		\[(1-r,1-t)H_L(p,q)(1-r,1-t)^t\ge c[(1-r)^2+2(1-r)(1-t)+(1-t)^2]=c(2-r-t)^2.\]
		Combining this with the observation that
		\[\lim_{r,t\to 1}\frac{(1-r-t)^2}{(2-r^2-t^2)^2}=\frac{1}{4}\]
		then completes the proof.
	\end{proof}
	
	\begin{remark}
		\label{remark:non-degen-new}
		By Lemma \ref{lem:upper bound on H_1,H^1 covariance}, (\ref{eqn:off-1}), (\ref{eqn:off-2}), and Corollary \ref{corr:ind cor} we see that for $r,t\in [-1,1]$, such that either $|r|<1$ or $|t|<1$, we have that $(\D h_N(\n),\D h_N(\n(r,t)))$ is a non-degenerate Gaussian random vector. As the covariance function for $h_N$ is isotropic, we see by applying a rotation, that as long as either $N^{-1}_1|(\sigma,\sigma')|<1$ or $N^{-1}_2|(\tau,\tau')|<1$, the law of
		\[(\D h_N(\sigma,\tau),\D h_N(\sigma',\tau'))\]
		is non-degenerate.
	\end{remark}
	
	We will also need a result that controls the rate at which the entries of the Hessian vanish as $r,t\to 1$.
	
	\begin{lemma}\label{lem:Covariance bounds around r,t=1}
		There is a fixed constant $C>0$, only dependent on $p,q$, such that for $1\le i\le N_1$, we have that
		\[\E_{\D h_N(\n),\D h_N(\n(r,t))}\bigg[\left(E_i E_1h_N(\n)\sq{1-r^2}+E_i E^1h_N(\n)\sq{1-t^2}\right)^2\bigg]\le C(2-t^2-r^2)^2,\label{Covariance bounds around r,t=1, epn 1}\]
		and for $1\le i\le N_2$, we have that
		\[\E_{\D h_N(\n),\D h_N(\n(r,t))}\bigg[\left(E_i E^1h_N(\n)\sq{1-r^2}+E^i E^1h_N(\n)\sq{1-t^2}\right)^2\bigg]\le C(2-t^2-r^2)^2\label{Covariance bounds around r,t=1, epn 2}\] 
	\end{lemma}
	\begin{proof}
		We only show (\ref{Covariance bounds around r,t=1, epn 1}), the case of (\ref{Covariance bounds around r,t=1, epn 2}) being similar. We also note that by symmetry and continuity, we may assume that $r,t\in (0,1)$. As the quantity on the right is positive, we may additionally assume that $r,t\in (1-\epsilon,1)$ for any fixed $\epsilon>0$.
		
		Let $\bar{h}_1$ and $\bar{h}_2$ be as in Lemma \ref{lem:hard-covariance-computation}, and let us denote $v_{N}(s)=(\sq{1-s^2},0,\dots 0)\in \R^{N-1}$ and $v(r,t)=(v_{N_1}(r),v_{N_2}(t))$. As both $\D \bar{h}_2(0)$ and $\D \bar{h}_1(v(r,t))$ are coordinate expressions of the gradient with respect to (potentially different) orthonormal frames, they are related at a fixed point by an invertible matrix. In particular, see that the linear relation $\D h_N(\n(r,t))=\D \bar{h}_2(0)=0$ is equivalent to the linear relation $\D \bar{h}_1(v(r,t))=0$. If we denote $\bar{h}_1=\bar{g}$ we see then that
		\[\E_{\D h_N(\n),\D h_N(\n(r,t))}[f]=\E[f|\D \bar{g}(0)=\D \bar{g}(v(r,t))=0].\label{eqn:abe-bird}\]
		With this notation, we see from the proof of Lemma \ref{lem:hard-covariance-computation} that we may rewrite (\ref{Covariance bounds around r,t=1, epn 1}) as
		\[\E_{\D h_N(\n),\D h_N(\n(r,t))}\left[\left((\sq{1-r^2},\sq{1-t^2})^t(\frac{d^2}{dx_1dx_i}\bar{g}(0),\frac{d^2}{dy_1dx_i}\bar{g}(0))\right)^2\right].\label{eqn:abe-ignore}\]
		We will write $(x,y)\in \R^{N_1-1}\times \R^{N_2-1}$ to denote the standard Euclidean coordinates. Let us denote, for $1\le i\le N_1$ \[H_{r,t}^i=\begin{bmatrix}\frac{d^3}{dx_1^2dx_i}\bar{g}(v(r,t))&\frac{d^3}{dx_1dy_1dx_i}\bar{g}(v(r,t))\\
			\frac{d^3}{dx_1dy_1dx_i}\bar{g}(v(r,t))&\frac{d^3}{dy_1^2dx_i}\bar{g}(v(r,t))\\
		\end{bmatrix}.\]
		We observe that by Taylor's Theorem, we have that
		\[|\frac{d}{dx_i}\bar{g}(v(r,t))-\frac{d}{dx_i}\bar{g}(0)-\frac{d^2}{dx_idx_1}\bar{g}(0)\sq{1-t^2}-\frac{d^2}{dx_idy_1}\bar{g}(0)\sq{1-r^2}|\le \]\[(\sup_{u,v\in [-1,1]^2}\|H_{u,v}^i\|)\|(\sq{1-r^2},\sq{1-t^2})\|^2.\]
		Observing that $\|(\sq{1-r^2},\sq{1-t^2})\|^2=(2-r^2-t^2)$ we thus see by (\ref{eqn:abe-bird}) that this implies that
		\begin{gather}
			\E_{\D h_N(\n),\D h_N(\n(r,t))}\left[\left((\sq{1-r^2},\sq{1-t^2})^t(\frac{d^2}{dx_1dx_i}\bar{g}(0),\frac{d^2}{dy_1dx_i}\bar{g}(0))\right)^2\right]\le \\
			(2-r^2-t^2)^2\E_{\D h_N(\n),\D h_N(\n(r,t))}[\sup_{u,v\in [-1,1]^2}\|H_{u,v}^i\|^2].
		\end{gather}
		Employing that $\|H_{u,v}^i\|\le \|H_{u,v}^i\|_{HS}$, and that conditioning on a Gaussian vector may only decrease its variance, we see that
		\[\E_{\D h_N(\n),\D h_N(\n(r,t))}[\sup_{u,v\in [-1,1]^2}\|H_{u,v}^i\|^2]\le \E[\sup_{u,v\in [-1,1]^2}\|H_{u,v}^i\|_{HS}^2].\label{RHS-Supremum}\]
		As the covariance of the entries of $H_{u,v}^i$ are clearly bounded (uniformly in $u,v$ and $N$), we see that the right-hand side of (\ref{RHS-Supremum}) is bounded by a fixed constant by the Borell-TIS inequality (see Theorem 2.1.1 of \cite{AT}). In view of (\ref{eqn:abe-ignore}), this completes the proof.
	\end{proof}
	
	Finally, we will end this section by computing a useful expression involving $\Sigma_U(r,t)$ that we will use above. We will fix $r,t\in (-1,1)$ for the remainder of the section. We first note that by Corollary \ref{corr:ind cor} and (\ref{appendix Sigma_E}) that
	\begin{gather}
		\Sigma_U(r,t)=\E[(h_N(\n),h_N(\n(r,t)))(h_N(\n),h_N(\n(r,t)))^t|E_1h_N(\n(r,t))=\\
		E^1h_N(\n(r,t))=E_1h_N(\n)=E^1h_N(\n)=0].\label{appendix-Sigma_E-2}
	\end{gather}
	To simplify this further, we observe that by direct computation from Lemma \ref{lem:hard-covariance-computation} the vector
	\[\frac{1}{\sq{2}}\left(h_N(\n)+h_N(\n(r,t)),E_1 h_N(\n)-E_1 h_N(\n(r,t)),E^1 h_N(\n)-E^1 h_N(\n(r,t))\right)\label{eqn:ignore-bing}\]
	and the vector
	\[\frac{1}{\sq{2}}\left(h_N(\n)-h_N(\n(r,t)),E_1 h_N(\n)+E_1 h_N(\n(r,t)),E^1 h_N(\n)+E^1 h_N(\n(r,t))\right)\]
	are independent. The covariance matrix of the vector (\ref{eqn:ignore-bing}) is given by
	\[\Sigma_Y(r,t):=\begin{bmatrix}
		1+r^pt^q& -pr^{p-1}t^qr_* & -q r^p t^{q-1}t_*\\
		-pr^{p-1}t^qr_* & p(1-r^pt^q+(p-1)r^{p-2}t^qr_*^2)& pq r^{p-1}t^{q-1}r_*t_*\\
		-q r^p t^{q-1}t_* & pq r^{p-1}t^{q-1}r_*t_*& q(1-r^pt^q+(q-1)r^{p}t^{q-2}t_*^2)
	\end{bmatrix}.
	\]
	We observe that by the above independence we see that
	\[-\frac{1}{2}(E,E)\Sigma_U(r,t)^{-1}(E,E)^t=-[\Sigma_Y(r,t)^{-1}]_{11}E^2.\]
	We now proceed to compute this quantity. We first note that by direct computation, we have that
	\[\det(\Sigma_Y(r,t))=pqb_{p,q}(r,t),\label{eqn:b-det-Y}\]
	\[\;b_{p,q}(r,t):=(1-r^{p}t^q)(1-r^{2p-2}t^{2q-2})+(p-1)(1-r^2)r^{p-2}t^{q}(1-r^{p}t^{q-2})\]\[+(q-1)(1-t^2)r^{p}t^{q-2}(1-r^{p-2}t^{q}).\]
	We will also need a term from the adjugate matrix, namely
	\[[\Sigma_Y(r,t)]_{22}[\Sigma_Y(r,t)]_{33}-[\Sigma_Y(r,t)]_{23}^2=pq b_{p,q}(r,t)-pq r^p t^{q}(1-r^{p-2}t^{q})(1-r^{p}t^{q-2}).\]
	Together these yield that
	\[[\Sigma_Y(r,t)^{-1}]_{11}=\frac{1}{\det(\Sigma_Y)}([\Sigma_Y(r,t)]_{22}[\Sigma_Y(r,t)]_{33}-[\Sigma_Y(r,t)]_{23}^2)=\]
	\[1-\frac{ r^p t^{q}(1-r^{p-2}t^{q})(1-r^{p}t^{q-2})}{b_{p,q}(r,t)},\]
	so altogether, we see that
	\[-\frac{1}{2}(1,1) \Sigma_U(r,t)^{-1}(1,1) ^t=-1+\frac{ r^p t^{q}(1-r^{p-2}t^{q})(1-r^{p}t^{q-2})}{b_{p,q}(r,t)}.\label{better U-def}\]
	
	\section{Proof of Lemma \ref{lem:KR-basic}\label{section:appendix-KR}}
	
	In this section, we provide a proof of Lemma \ref{lem:KR-basic}. The proof will be almost identical to that of Lemma 4 of \cite{pspin-second}, which may be immediately adapted to this case given if we prove the following result.
	
	\begin{lemma}\label{lem:Appendix-non-degeneracy}
		For $p,q\ge 5$, and $r,t\in (-1,1)$, we have that the Gaussian vector
		\[(h_N(\n),h_N(\n(r,t)),\D h_N(\n),\D h_N(\n(r,t)),\D^2h_N(\n),\D^2h_N(\n(r,t)))\label{vector-for-cov-appendix}\] is non-degenerate, up to degeneracies required by symmetries of the Hessian. That is, the vector is non-degenerate if one only takes elements of the Hessians above and on the main diagonal.
	\end{lemma}
	
	\begin{proof}[Proof of Lemma \ref{lem:KR-basic}]
		One may proceed almost identically to the proof of Lemma 4 of \cite{pspin-second}, replacing their Lemma 32 with our Lemma \ref{lem:Appendix-non-degeneracy}.
	\end{proof}
	
	The proof of Lemma \ref{lem:Appendix-non-degeneracy} will follow from two lemmas, both proven at the end of this section. To state the first of these, we first recall the (normalized) spherical $\ell$-spin glass model (see (1.1) and (4.1) of \cite{pspin-second}). For $\ell\ge 1$ and $N\ge 2$, the (normalized) $\ell$-spin glass model, is a smooth centered Gaussian random field $h_{N,\ell}:S^{N-1}\to \R$, with covariance
	\[\E[h_{N,\ell}(\sigma)h_{N,\ell}(\sigma')]=(\sigma,\sigma')^\ell.\]
	We see that this covariance is related to that of the $(p,q)$-spin glass model in that
	\[\E[h_{N_1,N_2,p,q}(\sigma,\tau)h_{N_1,N_2,p,q}(\sigma',\tau')]=\E[h_{N_1,p}(\sigma)h_{N_1,p}(\sigma')]\E[h_{N_2,q}(\tau)h_{N_2,q}(\tau')].\]
	For $r\in (-1,1)$, let us denote $\n(r)=\n_N(r)=(r,\sq{1-r^2},0,\dots,0)$ and $\n=\n(1)$. Furthermore, let $(E_i)_{i=1}^{N-1}$ denote the choice of (piece-wise) smooth orthonormal frame field on $S^{N-1}$ defined in Lemma 30 of \cite{pspin-second}, and define $(\D h_\ell,\D^2 h_\ell)$ as above.
	
	\begin{lemma}\label{lem:Appendix-p-spin-non-degeneracy}
		For $\ell\ge 5$, and $r\in (-1,1)$, the Gaussian vector
		\[(h_{N,\ell}(\n),h_{N,\ell}(\n(r)),\D h_{N,\ell}(\n),\D h_{N,\ell}(\n(r)),\D^2 h_{N,\ell}(\n),\D^2 h_{N,\ell}(\n(r)))\]
		is non-degenerate, up to degeneracies required by symmetries of the Hessian.
	\end{lemma}
	
	\begin{remark}
		\label{remark:non-degen}
		This is a similar to (and reliant on) Lemma 32 of \cite{pspin-second}, which establishes that for $\ell\ge 3$ and $r\in (-1,1)$, the Gaussian vector 
		\[(\D h_{N,\ell}(\n),\D h_{N,\ell}(\n(r)),\D^2 h_{N,\ell}(\n),\D^2 h_{N,\ell}(\n(r)))\]
		is non-degenerate, up to degeneracies required by symmetries of the Hessian.
	\end{remark}
	
	We observe that the $E$ defined in the proof of Lemma \ref{lem:hard-covariance-computation} is such that $(E_i)_{i=1}^{N_1-1}$ and $(E^i)_{i=N_1}^{N-2}$ are given by the extension to $S^{N_1-1}\times S^{N_2-1}$ of derivations acting only $S^{N_1-1}$ and $S^{N_2-1}$, respectively. Moreover, by construction, these are chosen to so that their restrictions coincide with the $E$ defined for above for Lemma \ref{lem:Appendix-p-spin-non-degeneracy} with $\ell=p$, $N=N_1$, and $r=r$ and with $\ell=q$, $N=N_2$ and $r=t$, respectively. In particular, we see that Lemma \ref{lem:Appendix-non-degeneracy} follows from Lemma \ref{lem:Appendix-p-spin-non-degeneracy} and Lemma \ref{lem:Appendix-non-degeneracy-general} below.
	
	\begin{lemma}\label{lem:Appendix-non-degeneracy-general}
		For $1\le i\le 2$, let $h_i$ denote a smooth centered Gaussian field defined on an open subset $U_i\subseteq \R^{n_i}$, with smooth covariance function $f_i$. Let $h$ denote the smooth centered Gaussian field defined on $U_1\times U_2$ with covariance function $f_1 f_2$. For some $\ell \ge 1$ and $1\le i\le 2$, let us choose some sequence of points $r_1^i,\dots r_\ell^i\in U_i$, and let us denote $r_{k,l}=(r_k^1,r_l^2)$. Then if we have, for both $l=1,2$, that the Gaussian vector comprised of entries
		\[(h_l(r_i^l),\D h_l(r_i^l),\D^2 h_l(r_i^l))_{i=1}^{\ell}\label{Appendix non degenerate vector 1}\]
		is non-degenerate up to degeneracies required by symmetries of the Hessian, then the Gaussian vector comprised of entries
		\[(h(r_{ij}),\D h(r_{ij}),\D^2 h(r_{ij}))_{1\le i,j\le \ell}\label{Appendix non degenerate vector 2}\]
		is non-degenerate up to degeneracies required by symmetries of the Hessian.
	\end{lemma}
	
	We now complete this section by giving the proofs of Lemmas \ref{lem:Appendix-p-spin-non-degeneracy} and \ref{lem:Appendix-non-degeneracy-general}.
	
	\begin{proof}[Proof of Lemma \ref{lem:Appendix-p-spin-non-degeneracy}]
		For the duration of this proof we will denote $h_N:=h_{N,\ell}$. To establish the desired claim, it is sufficient to show both that the vector 
		\[(h_{N}(\n),h_{N}(\n(r)),\D h_{N}(\n),\D h_{N}(\n(r)))\label{eqn:throw-deg}\] 
		is non-degenerate, and that the law of $(\D^2 h_{N}(\n), \D^2 h_{N}(\n(r)))$ conditioned on the event \[(h_{N}(\n)=h_{N}(\n(r))=0,\D h_{N}(\n)=\D h_{N}(\n(r))=0)\label{eqn:event-throw}\] 
		is non-degenerate. 
		
		For the ease of the reader, we will begin by recalling all the results of \cite{pspin-second} that we will use. For clarity, their $(f,p)$ coincides with our $(h_N,\ell)$. In their Lemma 12 they show that the law of $(\D h_{N}(\n),\D h_{N}(\n(r)))$ is non-degenerate, and that the covariance matrix of $(h_{N}(\n),h_{N}(\n(r)))$ conditional on the event $(\D h_{N}(\n)=\D h_{N}(\n(r))=0)$ is given by a matrix $\Sigma_U(r)$, which they show is invertible in Remark 31. Together these show that (\ref{eqn:throw-deg}) is non-degenerate. Now conditional on (\ref{eqn:event-throw}), their Lemma 13 shows that the only non-trivial correlations between entries of $(\D^2h_N (\n),\D^2h_N (\n(r)))$ are between $(E_{i} E_j h_N(\n),E_{i} E_j h_N(\n(r)))$, and those enforced by symmetry of the Hessian. Thus we are reduced to showing that the conditional law of $(E_{i} E_j h_N(\n),E_{i} E_j h_N(\n(r)))$ is non-degenerate for each $1\le i\le j\le N-1$. The case of $i,j<N-1$ is clear from their description in item (2) of Lemma 13, and the case of $i<N-1$ and $j=N-1$ follows from their description in item (3) and their proof of Lemma 32, where they show that $\Sigma_Z(r)$ is invertible. The remaining case is $i=j=N-1$, where the conditional covariance matrix is proportional to the matrix $\Sigma_Q(r)$ they define in (10.2). 
		
		In principle, one should be able to derive the invertibility of $\Sigma_Q(r)$ directly from this expression, but due to the technicality of this expression, we instead employ a more abstract argument. To perform this, we need only observe that the expression of $\Sigma_Q(r)$ is $N$-independent, so that to show that $\Sigma_Q(r)$ is invertible, it is thus sufficient to show that the conditional law of $(E_{N-1}E_{N-1}h_N(\n),E_{N-1}E_{N-1}h_N(\n(r)))$ is non-degenerate in the case of $N=2$. Henceforth, we assume that $N=2$ and omit $N$ from the notation.
		
		We observe that $h_\ell$ coincides in law with the Gaussian field defined for $(x,y)\in S^1$ by \[f(x,y)=\sum_{i=0}^{\ell}a_ix^iy^{\ell-i},\label{eqn:stop-ignore}\]
		where here $(a_i)_{i=0}^{\ell}$ are independent centered Gaussian random variables with $\E[a_i^2]=\binom{\ell}{i}$. We observe that by continuity of the derivatives of $f$ as a function of $(a_i)_{i=1}^\ell$, if the matrix $\Sigma_Q(r)$ is degenerate, any (deterministic) homogeneous polynomial of degree-$\ell$ on $S^1$ satisfying the linear relation (\ref{eqn:event-throw}) must satisfies a fixed non-trivial linear relationship between the elements of $(E_{1} E_1 f(\n),E_1 E_1 f(\n(r)))$. On the other hand, the polynomials
		\[f_{1}(x,y)=x^2(y\sq{1-r^2}-x r)^3,\;\;f_{2}(x,y)=x^3(y\sq{1-r^2}-x r)^2\]
		satisfy (\ref{eqn:throw-deg}) but in addition have $E_1 E_1 f_2(\n)=E_1 E_1 f_1(\n(r))=0$ and $E_1 E_1 f_1(\n)=E_1 E_1 f_2(\n(r))=2(1-r^2)^{3/2}$. Thus, no such relationship is possible, which shows that the law of $(E_{N-1}E_{N-1}h_N(\n),E_{N-1}E_{N-1}h_N(\n(r)))$ is non-degenerate when $N=2$ and thus in general.
	\end{proof}
	
	\begin{proof}[Proof of Lemma \ref{lem:Appendix-non-degeneracy-general}]
		This proof will consist of showing that the covariance matrix of (\ref{Appendix non degenerate vector 2}), after removing duplicate entries of the Hessians, is a submatrix of the Kronecker product of the covariance matrices of (\ref{Appendix non degenerate vector 1}), where again we have removed duplicate entries. Indeed, as the Kronecker product of two strictly positive definite matrices is strictly positive definite, the lemma then follows from the interlacing property for the eigenvalues of a submatrix of a symmetric matrix. 
		
		To demonstrate that one is a submatrix of the other, we will make repeated use of the identity (\ref{covariance of derivative formula}) above, which expresses the covariance between derivatives of a Gaussian field in terms of the derivatives of the covariance function. For ease of notation, we introduce the vector-valued operator $\D^{U}$, which is defined for a function $f$ in neighborhood of $x\in \R^N$ as
		\[\D^{U}f(x)=\left(\frac{d^2}{dx_idx_j}f(x)\right)_{1\le i\le j\le N},\]
		considered as an $N(N+1)/2$-dimensional vector with the lexicographic ordering on $(i,j)$. Then for $l=1,2$, we see from (\ref{covariance of derivative formula}) that covariance of the vector $(h_l,\D h_l,\D^{U} h_l)$ evaluated at two points $x,y\in \R^{n_l}$ is given by
		\[
		\Sigma^{l}(x,y):=\begin{bmatrix}
			f_l(x,y)& \D_1 f_l(x,y) & \D^{U}_1 f_l(x,y)\\
			\D_2 f_l(x,y)& \D_2 \D_1 f_l(x,y) & \D_2 \D^{U}_1 f_l(x,y)\\
			\D^{U}_2 f_l(x,y)& \D^{U}_2 \D_1 f_l(x,y) & \D^{U}_2 \D^{U}_1 f_l(x,y)\\
		\end{bmatrix},
		\label{eqn:hard-corr-block}
		\]
		where here $\D_1 f_l$ and $\D_{2} f_l$ denote the gradient of $f_l:U_l\times U_l\to \R$ taken with respect to the first and second copy of $U_l$, respectively (and similarly for $\D^{U}$) and where the matrix indices such that the index of $\D_2$ and $\D^{U}_2$ give the row number and the index of $\D_1$ and $\D^{U}_1$ give the column number. With this notation we see that the covariance of the matrix of (\ref{Appendix non degenerate vector 1}), with repeated entries removed, is up to reordering given by
		\[\Sigma^l:=\begin{bmatrix}\Sigma^l(r^l_i,r^l_j) \end{bmatrix}_{1\le i,j\le \ell}.\]
		To specify the covariance matrix of (\ref{Appendix non degenerate vector 1}), we note that covariance function of $h$ is given by $f_1f_2:U_1\times U_2\times U_1\times U_2\to \R$, or explicitly, for $(x,y,z,w)\in U_1\times U_2\times U_1\times U_2$ we have that
		\[\E[h(x,y)h(z,w)]=f_1(x,z)f_2(y,w).\]
		We are interested in the covariance matrix of the vector $(h,\D h,\D^{U} h)$ at two points $(x,y),$
		$(z,w)\in U_1\times U_2$. We will denote by $(\D_I,\D_{II})$ the decomposition of the gradient on $U_1\times U_2$ into the first and second factor, and similarly for $(\D_{I}^{2,U},\D_{II}^{2,U})$. Using this notation, we will also reorder and write the vector $(h,\D h,\D^{U} h)$ as $(h,\D_I h,\D^{U}_I h,\D_{II} h,\D_{I} \D_{II} h,\D^{U}_{II} h)$. With respect to this decomposition, we may write the covariance matrix in terms of the induced $6$-by-$6$ block decomposition as
		\[
		\Sigma(x,y,z,w):=\label{eqn:sigma-hard-def}\]
		\[\begin{bmatrix}
			
			f_1f_2& \D_1 f_1f_2 & \D^{U}_1 f_1f_2&
			f_1\D_1 f_2 & \D_1 f_1\D_1 f_2 & f_1\D^U_1 f_2\\
			
			\D_2 f_1f_2& \D_{12} f_1f_2 & \D^{U}_1 \D_2 f_1f_2&
			\D_2f_1\D_1 f_2 & \D_{12} f_1\D_1 f_2 & \D_2 f_1\D^U_1 f_2\\
			
			\D_2^U f_1f_2& \D_1\D_{2}^U f_1f_2 & \D^{U}_1 \D_2^U f_1f_2&
			\D_2^U f_1\D_1 f_2 & \D_1\D_{2}^U f_1\D_1 f_2 & \D_2^U f_1\D^U_1 f_2\\
			
			f_1\D_2f_2& \D_1 f_1\D_2f_2 & \D^{U}_1 f_1\D_2f_2&
			f_1\D_{12} f_2 & \D_1 f_1\D_{12} f_2 & f_1\D^U_1\D_2 f_2\\
			
			\D_2 f_1\D_2f_2& \D_{12} f_1\D_2f_2 & \D^{U}_1 \D_2 f_1\D_2f_2&
			\D_2 f_1\D_{12} f_2 & \D_{12} f_1\D_{12} f_2 & \D_2 f_1\D^U_1\D_2 f_2\\
			
			f_1\D_2^U f_2& \D_1 f_1\D_2^U f_2 & \D^{U}_1 f_1\D_2^U f_2&
			f_1\D_{1}\D_2^U f_2 & \D_1 f_1\D_{2}^U\D_1 f_2 & f_1\D^U_1\D_2^U f_2\\
			
		\end{bmatrix},
		\]
		where here $f_1=f_1(x,z)$, $f_2=f_2(y,w)$ (and similarly for their derivatives), $\D_{12}=\D_1\D_2$ and $\D_i f_l$ is given by its expression above. The covariance matrix of (\ref{Appendix non degenerate vector 2}), with degenerate entries removed, is then given by up to reordering by
		\[\Sigma:=\begin{bmatrix}\Sigma^l(r_i,r_j) \end{bmatrix}_{1\le i,j\le \ell}.\]
		Now the matrix $\Sigma^1(x,z)\otimes \Sigma^2(y,w)$, where $\otimes$ denotes the Kronecker product, may be decomposed into blocks labeled by $(i,j)$ with $1\le i,j\le 3$, by employing the block decomposition of (\ref{eqn:hard-corr-block}). We observe then that $\Sigma(x,y,z,w)$ is equivalent to the submatrix of $\Sigma_1(x,z)\otimes \Sigma_2(y,w)$, consisting of blocks with indices $(i,j)$ with $i+j\le 3$. As the choice of indices does not depend on the choice of $(x,y,z,w)$, we see as well that $\Sigma$ is a submatrix of $\Sigma^1\otimes \Sigma^2$.
	\end{proof}
	
	\section{Proof of Equation \ref{Single point complexity} and \ref{Single point complexity minima}\label{appendix:mckenna rescaling}}
	
	In this section, we will show that formulas (\ref{Single point complexity}) and (\ref{Single point complexity minima}) coincide with their counterparts in Theorem 2.1 of \cite{bipartite}. In particular, we need to relate our $\mu_{p,q,\gamma}$ to the $``\mu_{\infty}(u)''$ defined in Remark 2.2 and 2.3 of \cite{bipartite}, and in addition, we must relate the quantity $``E_{\infty}(p,q,\gamma)''$, defined in their Theorem 2.1 and Remark 2.3 to our $E_{p,q,\gamma;\infty}$. To avoid confusion, let us denote these quantities as $\mu_{p,q,\gamma,u}^{Mc}$ and $E_{p,q,\gamma;\infty}^{Mc}$, respectively. Specifically, to verify the desired claims of (\ref{Single point complexity}) and (\ref{Single point complexity minima}) we need to show that 
	\[E_{p,q,\gamma;\infty}=E_{p,q,\gamma;\infty}^{Mc};\;\;\int \log(|x|)\mu_{p,q,\gamma,u}^{Mc}(dx)=\int \log(|x-u|)\mu_{p,q,\gamma}(dx)+2C_{p,q,\gamma}-1,\label{appendix-eqn:mu}\]
	as well as verifying that the constructions given in Section \ref{subsection:notation} are indeed possible. In particular, given (\ref{appendix-eqn:mu}) we may reduce our definition of $\Sigma_{p,q,\gamma}$, given in (\ref{eqn:omega-def}) and (\ref{eqn:c-def}), to the complexity function of \cite{bipartite}.
	
	Before proceeding, we will explain the reason for the difference between our expressions and summarize the method to demonstrate (\ref{appendix-eqn:mu}). The difference is due to our rescaling of the rows and columns in the Hessian appearing in the Kac-Rice formula before passing to the limiting spectral measure. This trick, specific to the pure case, allows us to realize all the limiting empirical measures as translations of a single fixed measure. On the other hand, while the effect of our rescaling has an obvious effect on the determinant, it is less clear what effect it has on the limiting spectral measure, as there appears to be no direct way to compare the solution to our (\ref{system of equations}) to the solutions of (2.6) in \cite{bipartite}. Instead, we will use that this rescaling does not affect any of the proofs present in \cite{bipartite}, and so we obtain an alternative expression for the complexity, which, when combined with the expression in for this in \cite{bipartite} gives (\ref{appendix-eqn:mu}).
	
	To begin, we recall some results used in the proof of Theorem 2.1 in \cite{bipartite}. In particular, we need the following results, which follow from their application of the Kac-Rice formula in Lemma 3.2 (and noting that their $\alpha_1=\alpha_2=0$ in the pure case, so that $H_N(u)$ only depends on $u_0$): for nice $B\subseteq \R$,
	\[\lim_{N\to\infty}\frac{1}{N}\log(\E[\Crit_N(B)])=\frac{1}{2}[1+\gamma\log(\gamma/p)+(1-\gamma)\log((1-\gamma)/q)]+\]
	\[\lim_{N\to\infty}\frac{1}{N}\log(\int_B e^{-Nu^2/2}\E[|\det(H_N(u))|]du),\label{appendix-kac-rice-1}\]
	\[\lim_{N\to\infty}\frac{1}{N}\log(\E[\Crit_{N,0}(B)])=\frac{1}{2}[1+\gamma\log(\gamma/p)+(1-\gamma)\log((1-\gamma)/q)]+\]
	\[\lim_{N\to\infty}\frac{1}{N}\log(\int_B e^{-Nu^2/2}\E[|\det(H_N(u))|I(H_N(u)\ge 0)]du),\label{appendix-kac-rice-2}\]
	where here $H_N(u)$ is an $(N-2)$-by-$(N-2)$ Gaussian symmetric random matrix with independent entries, satisfying 
	\[\E[(H_N(u))_{ij}]=u\delta_{ij}(\frac{N}{N_1}p\delta_{i\le N_1-1}+\frac{N}{N_2}q\delta_{i>N_1-1})\]
	\[\Cov(H(u)_{ij})=
	\begin{cases}
		\frac{Np(p-1)}{N_1^2}(1+\delta_{ij});\;\;1\le i,j\le N_1-1\\
		\frac{Npq}{N_1N_2};\;\;1\le i\le N_1-1<j\le N-2\\
		\frac{Nq(q-1)}{N_2^2}(1+\delta_{ij});\;\; N_1-1< i,j\le N-2
	\end{cases}.
	\]
	Following their proof of Theorem 2.1 we see that
	\[\lim_{N\to\infty}\frac{1}{N}\log(\int_B e^{-Nu^2/2}\E[|\det(H_N(u))|]du)=\]
	\[\sup_{u\in B}\left(-\frac{u^2}{2}+\int \log(|x|)\mu_{p,q,\gamma,u}^{Mc}(dx)\right)\label{eqn:new-mc-1}\]
	\[\lim_{N\to\infty}\frac{1}{N}\log(\int_B e^{-Nu^2/2}\E[|\det(H_N(u))|I(H_N(u)\ge 0)]du)=\]
	\[\sup_{u\in (-\infty,-E_{p,q,\gamma;\infty}^{Mc})\cap B}\left(-\frac{u^2}{2}+\int \log(|x|)\mu_{p,q,\gamma,u}^{Mc}(dx)\right)\label{eqn:new-mc-2}\]
	Let us define the matrix \[\bar{D}_{N}=\begin{bmatrix}\sq{N_1/(Np)}I_{N_1-1}&0\\ 0&\sq{N_2/(Nq)}I_{N_2-1}\end{bmatrix},\] and define \[H^D_N(u)=\bar{D}_{N}H_N(u)\bar{D}_{N},\label{appendix:H^D_N def}\] 
	If we denote $H^D_N=H^D_N(0)$, we note that $H^D_N(u)$ coincides with $H_N^D-uI$ in law. We also note that 
	\[|\det(H_N^D-uI)|\stackrel{d}{=}\det(D_N)^2|\det(H_N(u))|=\exp(N(2C_{p,q,\gamma}-1+o(1)))|\det(H_N(u))|.\]
	Noting as well that the index is invariant under congruence transformations, we obtain the following result, which we formulate as a self-contained lemma for usage above.
	
	\begin{lemma}\label{lem:mckenna lem:G-det upperbound}
		For $p,q\ge 2$, $0<\gamma<1$, and nice $B\subseteq \R$ we have that
		\[\lim_{N\to\infty}\frac{1}{N}\log(\E[\Crit_N(B)])=C_{p,q,\gamma}+\lim_{N\to \infty}\frac{1}{N}\log(\int_B e^{-Nu^2/2}\E[|\det(H_N^D-uI)|]du),\]
		\[\lim_{N\to\infty}\frac{1}{N}\log(\E[\Crit_{N,0}(B)])=\]\[C_{p,q,\gamma}+\lim_{N\to \infty}\frac{1}{N}\log(\int_B e^{-Nu^2/2}\E[|\det(H_N^D-uI)|]I(H_N^D\ge uI)du),\]
		where $H^{D}$ is a symmetric random matrix with independent centered Gaussian entries with covariances for $1\le i,j\le N-2$ satisfying
		\[\E[(H_N^D)_{ij}^2]=\frac{1}{N}
		\begin{cases}
			p^{-1}(p-1)(1+\delta_{ij});\;\;1\le i,j\le N_1-1\\
			1;\;\;1\le i\le N_1-1<j\le N-2\\
			q^{-1}(q-1)(1+\delta_{ij});\;\; N_1-1< i,j\le N-2
		\end{cases}.
		\]
	\end{lemma}
	
	Now, proceeding through their Section 3, replacing $H_N(u)$ by $H_N^D(u)$, rescaling terms as necessary, the proofs in Section 3 of \cite{bipartite} give the following result.
	
	\begin{lemma}\label{lem:mckenna lem:G-det upper bound-first}
		Fix $p,q\ge 2$ and $0<\gamma<1$. Then there exists a unique compact-supported probability measure with continuous bounded density on $\R$, $\mu_{p,q,\gamma}$, specified by the procedure in (\ref{system of equations}). Moreover, for nice $B\subseteq \R$
		\[\lim_{N\to\infty}\frac{1}{N}\log(\int_B e^{-Nu^2/2}\E[|\det(H_N^D-uI)|]du)=\sup_{u\in B}\Theta_{p,q,\gamma}(u),\]
		\[\lim_{N\to\infty}\frac{1}{N}\log(\int_B e^{-Nu^2/2}\E[|\det(H_N^D-uI)|I(H_N^D-uI\ge 0)]du)=\]\[\sup_{u\in (-\infty,-E_{p,q,\gamma;\infty}^{Mc})\cap B}\Theta_{p,q,\gamma}(u).\]
	\end{lemma}
	
	Now employing this, and comparing the terms in Lemma \ref{lem:mckenna lem:G-det upperbound} and equations (\ref{eqn:new-mc-1}) and (\ref{eqn:new-mc-2}) we may conclude (\ref{appendix-eqn:mu}).


\begin{thebibliography}{10}
		
		\bibitem{AT}
		{\sc Adler, R.~J., and Taylor, J.~E.}
		\newblock {\em Random fields and geometry}.
		\newblock Springer Monographs in Mathematics. Springer, New York, 2007.
		
		\bibitem{bio4}
		{\sc Agliari, E., Barra, A., Bartolucci, S., Galluzzi, A., Guerra, F., and
			Moauro, F.}
		\newblock Parallel processing in immune networks.
		\newblock {\em Phys. Rev. E 87\/} (Apr 2013), 042701.
		
		\bibitem{cs2}
		{\sc Agliari, E., Barra, A., Galluzzi, A., Guerra, F., and Moauro, F.}
		\newblock Multitasking associative networks.
		\newblock {\em Phys. Rev. Lett. 109\/} (Dec 2012), 268101.
		
		\bibitem{MDE1}
		{\sc Ajanki, O.~H., Erd{\H{o}}s, L., and Kr{\"{u}}ger, T.}
		\newblock Stability of the matrix {D}yson equation and random matrices with
		correlations.
		\newblock {\em Probab. Theory Related Fields 173}, 1-2 (2019), 293--373.
		
		\bibitem{bio3}
		{\sc Amit, D.}
		\newblock {\em Modeling Brain Function: The World of Attractor Neural
			Networks}.
		\newblock Cambridge University Press, 1992.
		
		\bibitem{cupbook}
		{\sc Anderson, G.~W., Guionnet, A., and Zeitouni, O.}
		\newblock {\em An introduction to random matrices}, vol.~118 of {\em Cambridge
			Studies in Advanced Mathematics}.
		\newblock Cambridge University Press, Cambridge, 2010.
		
		\bibitem{exponential}
		{\sc Arous, G.~B., Bourgade, P., and McKenna, B.}
		\newblock Exponential growth of random determinants beyond invariance.
		\newblock {\em arXiv:2105.05000\/} (2021).
		
		\bibitem{subag-1RSB}
		{\sc Arous, G.~B., Subag, E., and Zeitouni, O.}
		\newblock Geometry and temperature chaos in mixed spherical spin glasses at low
		temperature: the perturbative regime.
		\newblock {\em Comm. Pure Appl. Math. 73}, 8 (2020), 1732--1828.
		
		\bibitem{pspin-one}
		{\sc Auffinger, A., Ben~Arous, G., and \v{C}ern\'{y}, J.}
		\newblock Random matrices and complexity of spin glasses.
		\newblock {\em Comm. Pure Appl. Math. 66}, 2 (2013), 165--201.
		
		\bibitem{tuca}
		{\sc Auffinger, A., and Chen, W.-K.}
		\newblock Free energy and complexity of spherical bipartite models.
		\newblock {\em J. Stat. Phys. 157}, 1 (2014), 40--59.
		
		\bibitem{azais}
		{\sc Aza\"{\i}s, J.-M., and Wschebor, M.}
		\newblock {\em Level sets and extrema of random processes and fields}.
		\newblock John Wiley \& Sons, Inc., Hoboken, NJ, 2009.
		
		\bibitem{11model}
		{\sc Baik, J., and Lee, J.~O.}
		\newblock Free energy of bipartite spherical {S}herrington-{K}irkpatrick model.
		\newblock {\em Ann. Inst. Henri Poincar\'{e} Probab. Stat. 56}, 4 (2020),
		2897--2934.
		
		\bibitem{bio1}
		{\sc {Barra}, A., and {Agliari}, E.}
		\newblock {A statistical mechanics approach to autopoietic immune networks}.
		\newblock {\em Journal of Statistical Mechanics: Theory and Experiment 2010}, 7
		(July 2010), 07004.
		
		\bibitem{socio2}
		{\sc Barra, A., and Contucci, P.}
		\newblock Toward a quantitative approach to migrants integration.
		\newblock {\em {EPL} (Europhysics Letters) 89}, 6 (Mar 2010), 68001.
		
		\bibitem{cs1}
		{\sc Barra, A., Genovese, G., Guerra, F., and Tantari, D.}
		\newblock How glassy are neural networks?
		\newblock {\em Journal of Statistical Mechanics: Theory and Experiment 2012},
		07 (Jul 2012), P07009.
		
		\bibitem{cs3}
		{\sc Barra, A., Genovese, G., Sollich, P., and Tantari, D.}
		\newblock Phase diagram of restricted boltzmann machines and generalized
		hopfield networks with arbitrary priors.
		\newblock {\em Phys. Rev. E 97\/} (Feb 2018), 022310.
		
		\bibitem{second1}
		{\sc Cammarota, V., Marinucci, D., and Wigman, I.}
		\newblock On the distribution of the critical values of random spherical
		harmonics.
		\newblock {\em J. Geom. Anal. 26}, 4 (2016), 3252--3324.
		
		\bibitem{second2}
		{\sc Cammarota, V., and Wigman, I.}
		\newblock Fluctuations of the total number of critical points of random
		spherical harmonics.
		\newblock {\em Stochastic Process. Appl. 127}, 12 (2017), 3825--3869.
		
		\bibitem{Parisi-General-Spherical}
		{\sc Chen, W.-K.}
		\newblock The {A}izenman-{S}ims-{S}tarr scheme and {P}arisi formula for mixed
		{$p$}-spin spherical models.
		\newblock {\em Electron. J. Probab. 18\/} (2013), no. 94, 14.
		
		\bibitem{Varadhan}
		{\sc Dembo, A., and Zeitouni, O.}
		\newblock {\em Large deviations techniques and applications}, vol.~38 of {\em
			Stochastic Modelling and Applied Probability}.
		\newblock Springer-Verlag, Berlin, 2010.
		\newblock Corrected reprint of the second (1998) edition.
		
		\bibitem{socio1}
		{\sc {Durlauf}, S.~N.}
		\newblock {How Can Statistical Mechanics Contribute to Social Science?}
		\newblock {\em Proceedings of the National Academy of Science 96}, 19 (Sept.
		1999), 10582--10584.
		
		\bibitem{MDE2}
		{\sc Erd{\H{o}}s, L., Kr{\"{u}}ger, T., and Schr{\"{o}}der, D.}
		\newblock Random matrices with slow correlation decay.
		\newblock {\em Forum Math. Sigma 7\/} (2019), Paper No. e8, 89.
		
		\bibitem{Guerra}
		{\sc Guerra, F.}
		\newblock Broken replica symmetry bounds in the mean field spin glass model.
		\newblock {\em Comm. Math. Phys. 233}, 1 (2003), 1--12.
		
		\bibitem{Log-Sobolev}
		{\sc Ledoux, M.}
		\newblock Concentration of measure and logarithmic {S}obolev inequalities.
		\newblock In {\em S\'{e}minaire de {P}robabilit\'{e}s, {XXXIII}}, vol.~1709 of
		{\em Lecture Notes in Math.} Springer, Berlin, 1999, pp.~120--216.
		
		\bibitem{bipartite}
		{\sc McKenna, B.}
		\newblock Complexity of bipartite spherical spin glasses.
		\newblock {\em arXiv:2105.05043\/} (2021).
		
		\bibitem{nic1}
		{\sc Nicolaescu, L.~I.}
		\newblock Critical sets of random smooth functions on compact manifolds.
		\newblock {\em Asian J. Math. 19}, 3 (2015), 391--432.
		
		\bibitem{nic2}
		{\sc Nicolaescu, L.~I.}
		\newblock Critical points of multidimensional random {F}ourier series: variance
		estimates.
		\newblock {\em J. Math. Phys. 57}, 8 (2016), 083304, 42.
		
		\bibitem{Parisi-General-Ising}
		{\sc Panchenko, D.}
		\newblock The {P}arisi formula for mixed {$p$}-spin models.
		\newblock {\em Ann. Probab. 42}, 3 (2014), 946--958.
		
		\bibitem{Bipartite-PD}
		{\sc Panchenko, D.}
		\newblock The free energy in a multi-species {S}herrington-{K}irkpatrick model.
		\newblock {\em Ann. Probab. 43}, 6 (2015), 3494--3513.
		
		\bibitem{Parisi-Vector-Spin}
		{\sc Panchenko, D.}
		\newblock Free energy in the mixed {$p$}-spin models with vector spins.
		\newblock {\em Ann. Probab. 46}, 2 (2018), 865--896.
		
		\bibitem{Parisi-Potts-Spin}
		{\sc Panchenko, D.}
		\newblock Free energy in the {P}otts spin glass.
		\newblock {\em Ann. Probab. 46}, 2 (2018), 829--864.
		
		\bibitem{bio2}
		{\sc Parisi, G.}
		\newblock A simple model for the immune network.
		\newblock {\em Proceedings of the National Academy of Sciences of the United
			States of America 87 1\/} (1990), 429--33.
		
		\bibitem{chi-ref}
		{\sc Simon, M.~K.}
		\newblock {\em Probability Distributions Involving Gaussian Random Variables: A
			Handbook for Engineers, Scientists and Mathematicians}.
		\newblock Springer-Verlag, Berlin, Heidelberg, 2006.
		
		\bibitem{pspin-second}
		{\sc Subag, E.}
		\newblock The complexity of spherical {$p$}-spin models---a second moment
		approach.
		\newblock {\em Ann. Probab. 45}, 5 (2017), 3385--3450.
		
		\bibitem{Parisi-Even-Spherical}
		{\sc Talagrand, M.}
		\newblock Free energy of the spherical mean field model.
		\newblock {\em Probab. Theory Related Fields 134}, 3 (2006), 339--382.
		
		\bibitem{Parisi-Even-Ising}
		{\sc Talagrand, M.}
		\newblock The {P}arisi formula.
		\newblock {\em Ann. of Math. (2) 163}, 1 (2006), 221--263.
		
	\end{thebibliography}
\end{document}